\documentclass[twocolumn]{autart}    % Enable this line and disable the 
\usepackage{amsmath,amssymb}
\usepackage{bm}
\usepackage[numbers]{natbib}
\usepackage{subcaption}
\usepackage{color}
\usepackage{graphicx}          % Include this line if your 
\newcommand{\edit}[1]{{\color{black}#1}}
\newcommand{\revise}[1]{{\color{black}#1}}

\begin{document}

\begin{frontmatter}
%\runtitle{Insert a suggested running title}  % Running title for regular 
                                              % papers but only if the title  
                                              % is over 5 words. Running title 
                                              % is not shown in output.

%\title{Linear System Analysis and Optimal Control of Gas Networks} % Title, preferably not more 
                                                % than 10 words.
\title{Linear System Analysis and Optimal Control of Natural Gas Dynamics in Pipeline Networks} % Let's plan to keep the title of the original submission and preprint

% Putting this into an acknowledgements section after the conclusions
% \thanks[footnoteinfo]{The authors are grateful to Cody W. Allen for numerous helpful discussions and to E. Olga Skowronek for drawing the network diagrams in Figures \ref{fig:bc_conf} and \ref{fig:net_graph}. This study was supported by the U.S. Department of Energy's Advanced Grid Modeling (AGM) project ``Dynamical Modeling, Estimation, and Optimal Control of Electrical Grid-Natural Gas Transmission Systems''. Research conducted at Los Alamos National Laboratory is done under the auspices of the National Nuclear Security Administration of the U.S. Department of Energy under Contract No. 89233218CNA000001. S. Shivakumar was supported by the National Science Foundation grants No. 1739990 and 1935453. Report number: LA-UR-23-24798.}

\author[lanl,cnls,asu_math]{Luke S. Baker}\ead{lsbaker@lanl.gov},  
\author[asu_mech]{Sachin Shivakumar}\ead{sshivak8@asu.edu},              
\author[asu_math]{Dieter Armbruster}\ead{dieter@asu.edu},
\author[asu_math]{Rodrigo. B. Platte}\ead{rplatte@asu.edu},
\author[lanl]{Anatoly Zlotnik}\ead{azlotnik@lanl.gov}

\address[lanl]{Applied Mathematics \& Plasma Physics, Los Alamos National Laboratory, Los Alamos, NM, U.S.A}                                    
\address[cnls]{Center for Nonlinear Studies, Los Alamos National Laboratory, Los Alamos, NM, U.S.A}         
\address[asu_mech]{School for Engineering of Matter, Transport and Energy, Arizona State University, Tempe, AZ, U.S.A.}    
\address[asu_math]{School of Mathematical and Statistical Sciences, Arizona State University, Tempe, AZ, U.S.A.}

                             % keyword list or with the 
                                          % help of the Automatica 
                                          % keyword wizard

\begin{abstract}                          % Abstract of not more than 200 words.
% We investigate nonlinear and adaptive linear control systems for modeling compressor-actuated dynamics in natural gas pipeline networks. 
We design nonlinear and adaptive linear model-predictive control (MPC) techniques to minimize operational costs of compressor-actuated dynamics in natural gas pipeline networks. We establish stability of the local linear system and derive rigorous bounds on error between the nonlinear and linear system solutions. These bounds are used to quantify conditions under which the linear MPC can substitute the nonlinear MPC without significant loss of predictive accuracy. Furthermore, we prove and numerically verify that the computational cost of the linear MPC is orders of magnitude lower than that of solving the baseline optimal control problem. Numerical simulations are performed on nontrivial networks to demonstrate that the proposed MPC can effectively adapt to varying load conditions while maintaining nearly 95\% optimality.
\end{abstract}

\end{frontmatter}

%*********************************************************************************************************************************************************************
\section{Introduction}
%*********************************************************************************************************************************************************************

\edit{Transient flows of natural gas in transmission pipelines are becoming increasingly variable. This is largely a consequence of gas-fired generators used to offset load fluctuations caused by the intermittent and uncertain nature of renewable electricity generation. The goal of transient optimization of natural gas pipeline networks is to maintain system performance within specified pressure and flow constraints while minimizing the cost of compressor operation, in which pressure, flow rates, and compressor actions change over time in response to demand variations.  While steady-state optimization remains essential for capacity analysis and long-term planning \cite{wong1968optimization,percell1987steady,de2000gas,misra2014optimal}, optimization of transient flows is becoming critical for predictive day-ahead operations \cite{abbaspour2008nonisothermal,zlotnik2015optimal}. 
Indeed, the application of transient optimization for day-ahead planning offers an effective method for establishing compressor control policies and flow schedules that can serve as an operational baseline for the upcoming 24-hour period \cite{zlotnik2015optimal, zlotnik2015model}.

 Transient optimization of pipeline flows is typically formulated as a conventional nonlinear optimization problem, which is obtained by discretizing the cost objective and the governing partial differential equations (PDEs) of an underlying optimal control problem (OCP) over time and space \cite{ehrhardt2005nonlinear, steinbach2007pde, domschke2011combination, liu2011coordinated}. While these approaches have demonstrated the ability to enhance state estimation \cite{sundar2018state} and coordination with electric power grids \cite{zlotnik2016coordinated}, the collective optimization over all time samples at once may become problematic for computer memory and performance when applied to large-scale systems.  As an alternative, model-predictive control (MPC) has been proposed decades ago and continues to gain traction because of its potential to drastically reduce the number of optimization variables \cite{marques1988line, aalto2006real, gopalakrishnan2013economic, aalto2015model}. Unlike conventional transient optimization, MPC divides the planning horizon into a number of small time intervals and optimizes locally in time with a moving frame. Despite numerous advancements, the nonlinearity and poor conditioning of the optimization problems can still pose significant challenges for contemporary computational methods \cite{srinivasan2022numerical}.

 Linearization, model reduction, and preconditioning are powerful techniques that may be applied to improve the performance of numerical methods and address some of the challenges mentioned above \cite{grundel2013computing, behbahani2010accuracy, luongo1986efficient, feldman1988optimization, Yue20, srinivasan2022numerical}.  In some cases, application of such techniques can promote convergence of a numerical method to a feasible solution even when the original formulation fails \cite{srinivasan2022numerical}. Mixed-integer linear programming and sequential quadratic programming are generally very effective in generating physically meaningful solutions for the transient optimization problem \cite{domschke2011combination, gugat2018mip, beylin2020fast}. However, while pipeline field data has supported the application of linear system modeling for certain transient flows \cite{kralik1984modeling, allen2022gas}, recent simulations highlight that large variations in flow rates and amplitudes can undermine the validity of linear modeling \cite{baker2021analysis, baker2023gas}. Consequently, a rigorous and quantifiable error bound becomes indispensable for establishing conditions under which linear models remain accurate within a tolerable error.

This paper provides a comprehensive analysis of modeling and actuating gas flow with adaptive linear control systems.  First, transfer functions are used to quantify rates of variation for which simplifications of the governing PDEs can be made. The adaptive linear system is derived by linearizing a spatially discretized system about a nominal state. Crucially, the nominal state can be either steady or unsteady to enable adaptive re-linearization about unsteady instantaneous states in predictive controller applications.  A key contribution of this paper that shows the relevance of MPC to practical pipeline operations is rigorous bounding of the error between the solutions of the finite-dimensional linearized and nonlinear systems, expressed in terms of elapsed time and magnitude of variation about the nominal state. Because the result requires stability of the linear system, we prove that the linearized system is asymptotically stable and compare the eigenvalues of the state matrix to the poles of the transfer function. 

The second key contribution of this paper is the development of a versatile MPC synthesis method aimed at minimizing the energy expended by compressor units. As mentioned above, MPC typically divides the time horizon into multiple smaller segments and optimizes over neighboring future segments to determine the control policy for the current segment. Generally, the optimality of the solution improves as the number of neighboring segments increases. However, this comes with a trade-off between optimality and computational cost.  Rather than incorporating forecasts of future conditions and optimizing over a larger number of time segments, the significance of our proposed MPC formulation lies in its ability to make accurate predictions while optimizing only one time segment into the future. We prove and numerically verify reduction in computational costs by orders of magnitude compared to conventional optimal control (OC) formulations \cite{zlotnik2015optimal, baker2023optimal}. Numerical results for two nontrivial networks are provided to examine the trade-offs between MPC and OC in terms of optimality, computational time, and operational responsiveness.}

The remainder of the paper is outlined as follows.  The flow equations and transfer function representations are reviewed in Section \ref{sec:pipe_flow} for a single pipe and then extended to networks in Section \ref{sec:network}.  A lumped element discretization method is reviewed in Section \ref{sec:discretization} and an adaptive linear control system is derived in Section \ref{sec:linear_ti}.  Section \ref{sec:eigenvalues} establishes several results pertaining to the eigenvalues of the finite-dimensional state matrix and poles of the irrational transfer matrix.  Section \ref{sec:error} examines the error between the solutions of the nonlinear and linearized systems.  The MPC and OC formulations are presented in Section \ref{sec:MPC} and several computational studies are conducted in Section \ref{sec:results}. Finally, a review of compelling future research directions is highlighted in Section \ref{sec:conclusion}.

%*********************************************************************************************************************************************************************
\section{Pipeline Flow Equations}  \label{sec:pipe_flow}
%*********************************************************************************************************************************************************************

The flow of natural gas in a horizontal pipe is approximated with the one-dimensional, isothermal Euler equations given by \cite{osiadacz1984simulation,osiadacz2001comparison,thorley1987unsteady}
\begin{subequations} \label{eq:gaspde0}
\begin{align}
    \partial_t \rho + \partial_x (\rho v) & = 0, \label{eq:gaspde0a} \\
    \partial_t (\rho v) + \partial_x (\sigma^2 \rho + v^2 \rho) & = - \frac{\lambda}{2D}\rho v |v|.  \label{eq:gaspde0b} 
\end{align}
\end{subequations}
The variables $v(t,x)$ and $\rho(t,x)$ represent the velocity and density of the gas, respectively, at time $t\in [0,T]$ and axial location $x\in[0,\ell]$, where $T$ denotes the time horizon and $\ell$ denotes the length of the pipe.  The parameter $\sigma$ represents the sound speed that is generally dependent on the temperature and pressure of the gas in accordance to equations of state.  The above PDE system describes mass conservation \eqref{eq:gaspde0a} and momentum conservation \eqref{eq:gaspde0b}.   The Darcy-Weisbach term on the right-hand side of equation \eqref{eq:gaspde0b}  models momentum loss caused by turbulent friction in which $\lambda$ is called the friction factor.  Restrictions on initial and boundary conditions under which the model in \eqref{eq:gaspde0} is well-posed have been previously examined \cite{gugat2012well}.  Because analytical solutions of the isothermal Euler equations exist only for certain cases of initial and boundary conditions \cite{gugat2017isothermal}, the solution is obtained numerically for general systems and conditions.  Initial and boundary conditions for general network topologies are described below.  \edit{Throughout, we assume regular operating conditions for the intended application in which flow is subsonic and density is positive.}

\subsection{Simplified Equations}

The equation of state is approximated by the ideal gas law for which $\sigma$ in equation \eqref{eq:gaspde0b} is constant. Although non-ideal models may be needed for precise flow quantification in large gas pipelines, ideal gas modeling provides a qualitative description of the flow phenomena. We further suppose that the gas velocity $v$ is much lower than the sound speed $\sigma$ so that the nonlinear term $\partial_x (v^2\rho)$ may be removed from the momentum equation.  Moreover, the inertia term $\partial_t (\rho v)$ is omitted under the assumption of slowly-varying flows, resulting in the so-called friction-dominated PDE system \cite{gugat2020closed}.  These simplifications give rise to the flow equations
\begin{subequations} \label{eq:simplified}
\begin{align}
\partial _t \rho +\partial_x \varphi  &=0, \label{eq:simplified1} \\
\delta \partial_t \varphi +\sigma^2 \partial_x\rho &= -\frac{\lambda}{2D}\frac{ \varphi|  \varphi|}{\rho},  \quad \label{eq:simplified2}
\end{align}
\end{subequations}
where $\varphi=\rho v$ is the cross-sectional mass flux.  In the above equation, we use a binary parameter $\delta$ with $\delta = 1$ if the flux derivative term is included and $\delta = 0$ otherwise.  In the subsequent analysis, we consider both forms of simplified models and refer to the case in which $\delta =0$ as the friction-dominated approximation.

%*********************************************************************
\subsection{Transfer Matrix Analysis}  \label{sec:transfer_matrix}

Simulations suggest that the friction-dominated approximation is valid when boundary conditions vary sufficiently slowly over time \cite{osiadacz1984simulation, osiadacz2001comparison}.  Because results of numerical simulations do not generalize well to variations in system parameters, we propose an analytical method to quantify acceptable rates of variation in terms of pipeline parameters.  Assuming an initially steady-state system, we analyze the response to time-varying boundary conditions using frequency-domain transfer functions \cite{kralik1984modeling, zecchin}. These functions are derived by linearizing the simplified models about the steady state, approximating space-varying coefficients with nominal constants, and applying the Laplace transform to the resulting linear PDEs.

\begin{figure*}[!t]
    % \centering
    \begin{subfigure}[t]{0.5\textwidth}
        \centering
        \includegraphics[width=1\linewidth]{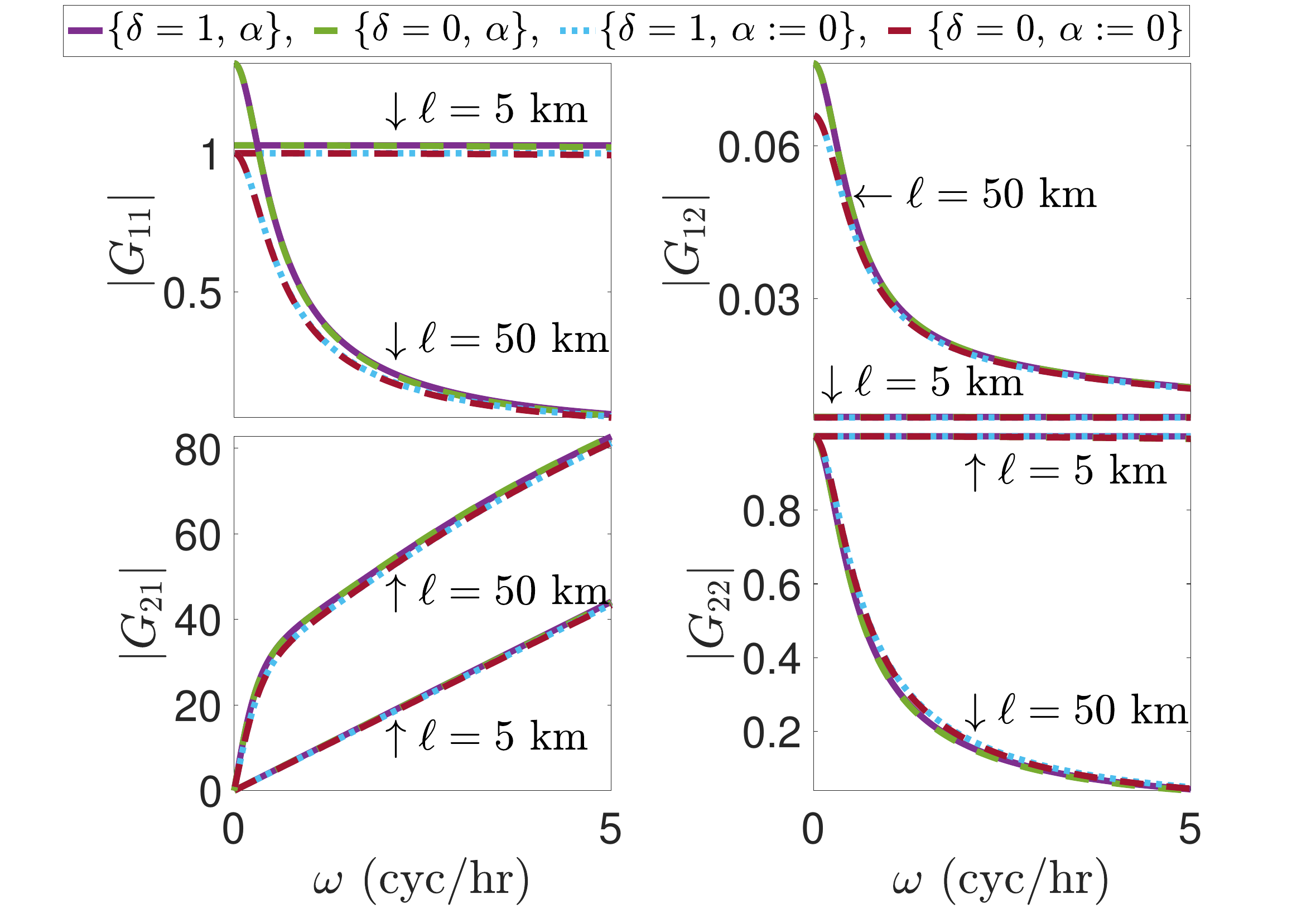}
        \caption{Magnitude}
    \end{subfigure}%
    ~ 
    \begin{subfigure}[t]{0.5\textwidth}
        \centering
        \includegraphics[width=1\linewidth]{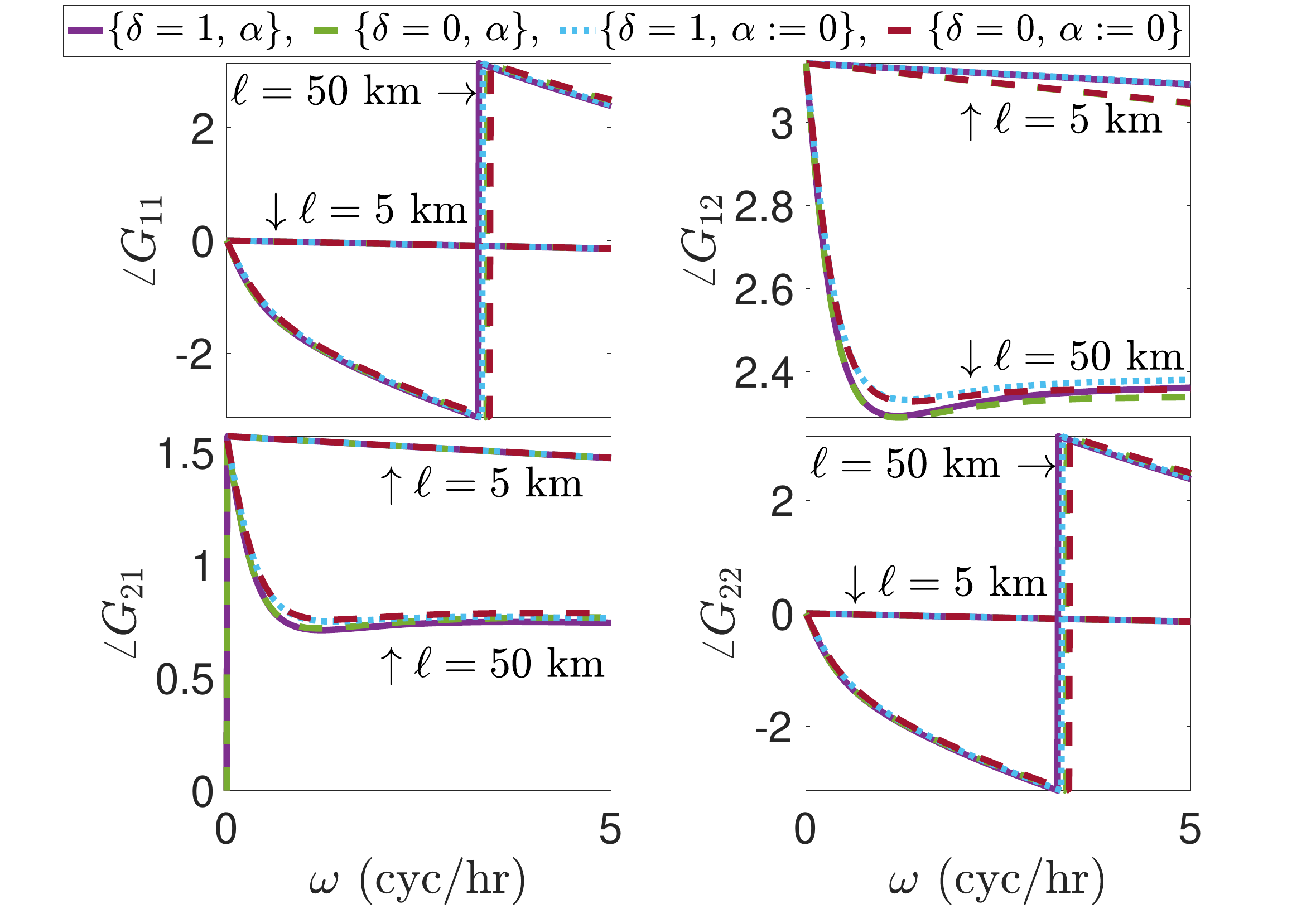}
        \caption{Phase}
    \end{subfigure}
    \caption{Magnitude and phase of the transfer matrix coefficients for two pipes of lengths 5 km and 50 km, and the four linear models produced by the combinations of taking zero and nonzero values for $\alpha$ and $\delta$, as indicated in each panel.  The remaining parameters are $D=0.5$ m, $\sigma=377$ m/s, and $\lambda=0.011$, with nominal values $\overline \varphi=300$ kg/sm$^2$ and $\overline \rho=35$ kg/m$^{3}$.}
    \label{fig:magG}
\end{figure*}

Extending our previous study \cite{baker2021analysis}, the above operations result in a representation of the form
\begin{subequations} \label{eq:linear_dyn}
\begin{align}
    s P +\partial_x \Phi  &=0, \label{eq:linear_dyn1} \\
\delta s\Phi +\sigma^2 \partial_xP &= \alpha P-\beta \Phi, \label{eq:linear_dyn2}
\end{align}
\end{subequations}
where $s$ is the complex frequency. Laplace transforms of variations of $\rho(t,x)$ and $\varphi(t,x)$ about nominal steady-state values as functions of $x$ are denoted by $P(s,x)$ and $\Phi(s,x)$, respectively. The true Jacobian coefficients are approximated with the constant parameters
\begin{equation} \label{eq:alpha_beta}
\alpha =\frac{\lambda \overline \varphi |\overline \varphi|}{2D \overline \rho}, \qquad 
\beta =\frac{\lambda |\overline \varphi|}{D \overline \rho},
\end{equation}
where $\overline \rho$ and $\overline \varphi$ are constant baseline values representing nominal density and steady-state mass flux, respectively.  Field data \cite{kralik1984universal} and simulations \cite{baker2021analysis} indicate that the term involving $\alpha$ in equation \eqref{eq:linear_dyn2} may be removed to further simplify the frequency domain representation without significant loss in accuracy.  Effects of this additional simplification will be analyzed below in Figure \ref{fig:magG}.  For ease of exposition, we denote inlet and outlet pipeline quantities with superscripts ``0" and ``$\ell$", respectively.   In analogy to an electric circuit, we define series impedance $z(s)=\delta s-\beta$ and propagation factors $\gamma_{\pm}(s)=\alpha/(2\sigma^2)\pm \gamma(s)$, where $\gamma(s)=\sqrt{\alpha^2/(4\sigma^4)+sz(s)/\sigma^2}$.  The solution of equations \eqref{eq:linear_dyn1}-\eqref{eq:linear_dyn2} evaluated at the ends of the pipe may be written in terms of these circuit quantities as \cite{baker2021analysis}
\begin{subequations} \label{eq:transfer_rep} 
\begin{align}
\!\!\!\!\! P^{\ell}&=  \frac{\gamma_{+} e^{\gamma_{+}\ell}-\gamma_{-} e^{\gamma_{-}\ell}}{\gamma_{+}-\gamma_{-}} P^0 -\frac{z}{\sigma^2}  \frac{ e^{\gamma_{+}\ell}- e^{\gamma_{-}\ell}}{\gamma_{+}-\gamma_{-}}  \Phi^0, \qquad \label{eq:transfer_rep_P} \\
\!\!\! \Phi^{\ell} &= -s \frac{ e^{\gamma_{+}\ell}- e^{\gamma_{-}\ell}}{\gamma_{+}-\gamma_{-}} P^0 -    \frac{\gamma_{-} e^{\gamma_{+}\ell}-\gamma_{+} e^{\gamma_{-}\ell}}{\gamma_{+}-\gamma_{-}}   \Phi^0. \label{eq:transfer_rep_Phi} 
\end{align}
\end{subequations}
Two of the four inlet and outlet variables are specified as boundary conditions, while the remaining two variables represent the responses at the opposite ends of the pipe.  For a single pipe, designated boundary conditions specify density and mass flux at opposite ends.  By rearranging the above relations as a mapping from input to output variables, the transfer matrix representation assumes the form
\begin{equation} \label{eq:transfer_matrix}
\begin{bmatrix}
P^{\ell} \\
\Phi^0
\end{bmatrix}
=
G
\begin{bmatrix}
P^0 \\
\Phi^{\ell}
\end{bmatrix},
\quad 
G(s)=\begin{bmatrix}
G_{11}(s) & G_{12}(s) \\
G_{21}(s) & G_{22}(s)
\end{bmatrix},
\end{equation}
where the coefficients of the $2\times 2$ transfer matrix $G(s)$ are uniquely determined by the coefficients in equation \eqref{eq:transfer_rep}.  Under the further assumption that $\alpha=0$, the transfer matrix simplifies to
\begin{equation}
    G(s)=
\begin{bmatrix}
\text{sech} \left(\ell \gamma(s) \right) & -\frac{z_c(s)}{\sigma} \text{tanh} \left(\ell \gamma(s)  \right) \\
\frac{\sigma}{z_c(s)}\text{tanh}\left(\ell \gamma(s)  \right) & \text{sech} \left(\ell \gamma(s)  \right)
\end{bmatrix},
\end{equation}
where $z_c(s)=\sqrt{z(s)/s}$.  

By evaluating equation \eqref{eq:transfer_matrix} for all combinations of zero and nonzero values for $\alpha$ and $\delta$, four linear systems are obtained. Formally, we suppose these systems are approximately equivalent if they attenuate and delay boundary condition variables in similar ways. To quantify suitable rates of variation, we analyze the magnitudes $|G_{mn}(\bm j \omega)|$ and phases $\angle G_{mn}(\bm j \omega)$ of the coefficients $G_{mn}(s)$ as functions of frequency $s=\bm j \omega$, where $\bm j=\sqrt{-1}$ is the imaginary unit and $\omega$ is the angular frequency. The magnitude describes amplitude attenuation and the phase introduces time delays between input and output waves. The systems are considered approximately equivalent over a range of frequencies if the magnitudes and phases of their transfer matrix coefficients are similar over this range. This equivalence can be rigorously quantified by setting a threshold on the absolute differences of magnitudes and phases between any two systems and inverting the differences numerically to express the maximum frequency of the allowable range as a function of the chosen threshold.

Figure \ref{fig:magG} displays the magnitudes and phases of the transfer matrix coefficients for each of the four linear systems as functions of input frequency. For short pipes, these coefficients closely align at frequencies below 5 cyc/hr. For longer pipes, the coefficients generally agree, except for $G_{11}(\bm j \omega)$ and $G_{21}(\bm j \omega)$, associated with the simplified system described by $\delta = 0$ and $\alpha = 0$. It is important to observe that the magnitude of $G_{21}(\bm j \omega)$ significantly increases with frequency.  This implies that pipeline systems could possibly be exposed to unsafe operating conditions if compressor control protocols are rapidly varying, which in the above analysis would correspond to high frequencies and sizable amplitudes in boundary conditions.  Such destabilization of the transfer function representation of gas pipeline flow was numerically analyzed in our previous study \cite{baker2021analysis}.  In practice, such variations are avoided by using only slowly-varying control actions, which makes the actions of $G_{11}$ and $G_{21}$ negligible.  Under these conditions, Figure \ref{fig:magG} demonstrates that all four linear models agree for mass flow withdrawal frequencies below 5 cyc/hr.  

\edit{The results for input density variations may be extended to infer appropriate operating conditions for compressor units.  As defined below, the action of a compressor is modeled as a multiplicative ratio between discharged and suctioned pressures at the inlet of the pipeline.  To maintain sufficiently low frequencies in the pressure discharged from a compressor unit,  either the compressor action and input density should both be composed of sufficiently low frequencies, or the compressor actions should be chosen to offset input density variations (e.g. completely out of phase).  In addition to possibly being infeasible, this second option conflicts with the desire to achieve minimal energy control. Therefore, we conclude that compressor actions and input density variations should both be composed of sufficiently low frequencies.}  We adhere to the frequency restrictions in our subsequent analysis.

%*********************************************************************************************************************************************************************
\section{Network Flow Equations}  \label{sec:network}
%*********************************************************************************************************************************************************************

We now extend the flow equations from a single pipe to a network of interconnected pipes. \edit{Although the MPC formulation defined below may be tailored to accommodate time-varying topological conditions, we assume that the network topology is fixed.}  A pipeline network is modeled as a connected and directed graph $(\mathcal E, \mathcal V)$ consisting of a set of edges $\mathcal E =\{1,\dots,E\}$ and a set of nodes $\mathcal V=\{1,2,\dots,V\}$, where $E$ and $V$ denote the numbers of edges and nodes of the graph. The set of edges represents the pipelines of the topology.  The node set is composed of junctions and stations where gas is injected into or withdrawn from the network. The elements of both sets are assumed to be ordered according to the integer labels, where the label $k$ is reserved for indexing edges in $\mathcal E$, and the labels $i$ and $j$ are reserved for indexing nodes in $\mathcal V$.   The nodes of the network are partitioned into nonempty sets of supply nodes $ \mathcal V_s\subset \mathcal V$ and withdrawal nodes $ \mathcal V_w\subset \mathcal V$. Supply and withdrawal nodes are ordered in $\mathcal V$ so that $i<j$ for all $i\in \mathcal V_s$ and $j\in \mathcal V_w$.  The graph is directed by assigning a positive flow direction along each edge with the convention that edges incident to supply nodes are directed away from them.  
The notation $k:i\mapsto j$ means that edge $k\in \mathcal E$ is directed from node $i\in \mathcal V$ to node $j\in \mathcal V$.  For each node $j\in \mathcal V$, we define (potentially empty) incoming and outgoing sets of pipelines by $_{\mapsto} j=\{k\in \mathcal E \mid k:i\mapsto j\}$ and $j_{\mapsto}=\{k\in \mathcal E \mid k:j\mapsto i \}$, respectively.  

For each pipe $k\in \mathcal E$, the state variables are cross-sectional density $\rho_k(t,x)$ and mass flux $\varphi_k(t,x)$ for time $t\in [0,T]$ and axial location $x\in[0,\ell_k]$, where $\ell_k$ denotes the length of the pipe.  As above, the flow of gas through the pipe indexed by $k\in \mathcal E$ is governed by
\begin{subequations} \label{eq:net_flow}
\begin{align}
\partial _t \rho_k +\partial_x \varphi_k  &=0, \label{1} \\
\delta\partial_t \varphi_k +\sigma^2 \partial_x\rho_k &= -\frac{\lambda_k}{2D_k}\frac{ \varphi_k|  \varphi_k|}{\rho_k}, \label{2} \quad
\end{align}
\end{subequations}
where $D_k$ and $\lambda_k$ are the diameter and friction factor of the pipe, respectively, and $\sigma$ is the sound speed.   As in the previous section, $\delta$ can be set to either zero or one, provided that the power spectrum of variations in the boundary conditions is negligible for frequencies above 5 cyc/hr, which is typically the case in the operating regime of gas transmission pipelines.

We define $\mathcal C \subset \mathcal E$ to be the set of edges $k\in \mathcal E$ incident to a compressor station and simply refer to an edge $k\in \mathcal C$ as a compressor.  By bisecting an edge into two edges if necessary, we may assume that each compressor $k\in \mathcal C$ is located at the inlet of its incident edge $k\in \mathcal E$ with respect to the graph orientation.  The action of the compressor $k\in \mathcal C$ is modeled with the multiplicative control variable $\mu_k(t)$ with $1\le \mu_k(t)$ for all $t\in [0,T]$.   In particular, the pressure of gas discharged from the compressor unit $k\in \mathcal C$ is $\mu_k(t)$ times larger than the suction pressure.  To simplify notation, we define $\mu_k=1$ for all $k\in \mathcal E \setminus \mathcal C$. 
\begin{figure}
\centering
\includegraphics[width=.8\linewidth]{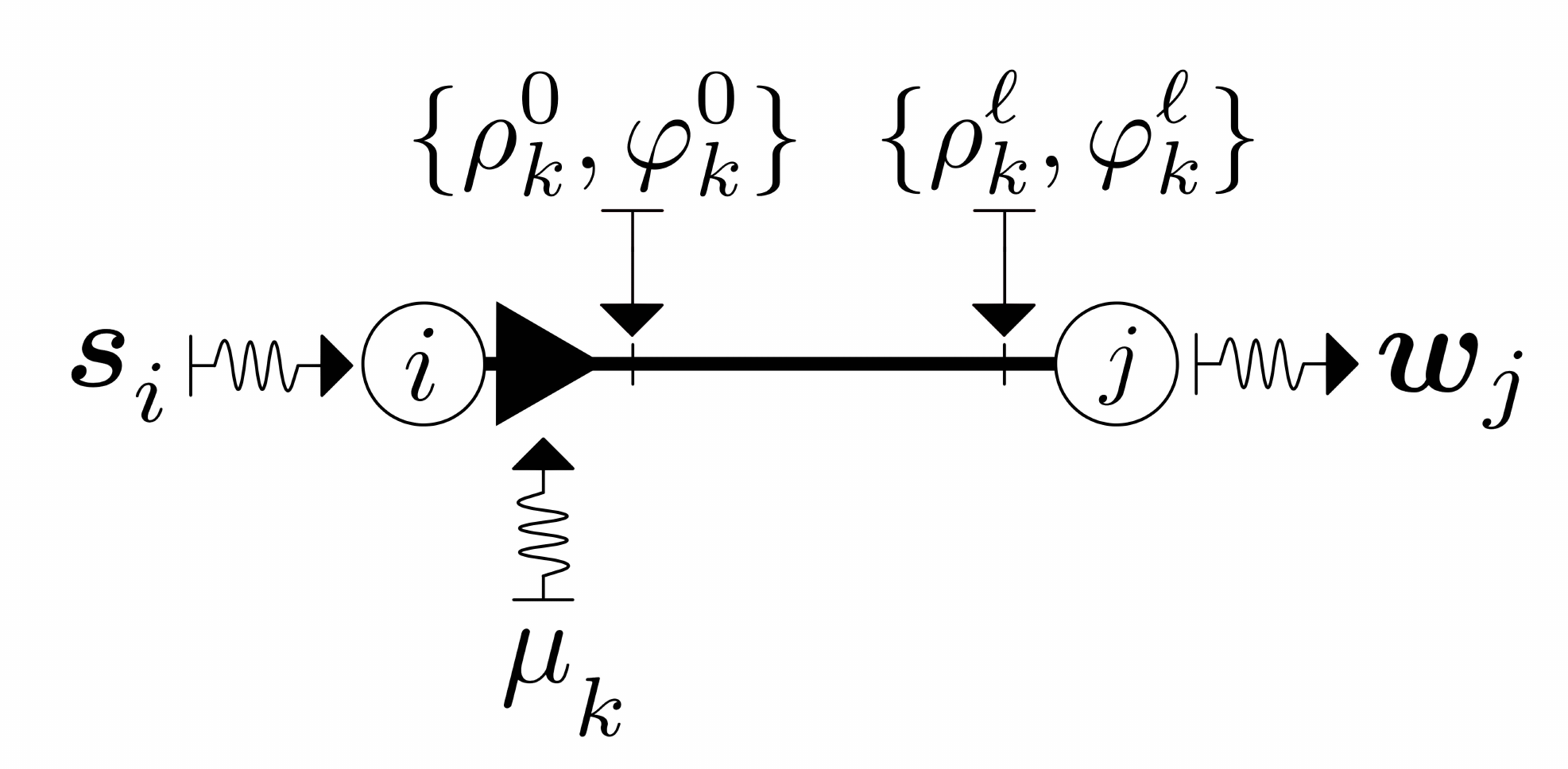}
\caption{Configuration of a pipeline segment $k:i\mapsto j$ with $i\in \mathcal V_s$ and $j\in \mathcal V_w$. Flow and control variables are indicated.}
\label{fig:bc_conf}
\end{figure}

Natural gas is injected into the network at supply nodes $i\in \mathcal V_s$, where supply density $\bm s_i(t)$ is specified as a time-varying boundary condition. While time-varying supply densities are considered in the derivations that follow, we emphasize that flows are typically simulated under the assumption of a constant supply density. The mass outflow from each withdrawal node $j\in \mathcal V_w$ is also defined as a boundary condition, represented by $\bm w_j(t)$. Throughout this study, all nodal quantities are denoted using bold symbols. In alignment with the input and output variables of the transfer matrix, inlet and outlet edge variables are denoted with superscripts ``0" and ``$\ell$", respectively, as indicated in Fig. \ref{fig:bc_conf}. Additionally, the cross-sectional area of pipe $k\in \mathcal E$ is defined by $\chi_k=\pi D_k^2/4$.  With these definitions, the boundary conditions governing gas flow in the network are expressed as
\begin{subequations} \label{eq:boundary_conditions}
\begin{align}
 \rho_{k}(t,0)&=\mu_{k}(t) \bm s_i(t),\quad 
\rho_{k}(t,\ell_k)= \bm \rho_{j}(t), \label{eq:bc1} \\
 \rho_{k}(t,0)&=\mu_{k}(t) \bm \rho_i(t),\quad 
\rho_{k}(t,\ell_k)= \bm \rho_{j}(t), \label{eq:bc2} \\
 \bm w_j(t)&=\sum_{k\in _{\mapsto}j} \chi_k \varphi_{k}^{\ell}(t)
 -\sum_{k\in j_{\mapsto}} \chi_k \varphi_{k}^0(t),  \label{eq:bc}
\end{align}
\end{subequations}
where equation \eqref{eq:bc1} is defined for  $k:i\mapsto j$ with $i\in \mathcal V_s$, equation \eqref{eq:bc2} is defined for $k:i\mapsto j$ with $i\in \mathcal V_w$, and equation \eqref{eq:bc} is defined for $j\in \mathcal V_w$.   These conditions incorporate the effects of compression and the conservation of mass flow through withdrawal nodes.  The initial condition for the network is defined by
\begin{equation}
 \rho_k(0,x)=\varrho_k(x), \quad \varphi_k(0,x)=\phi_k(x), \quad \forall k\in \mathcal{E}, \label{ic}
\end{equation}
\edit{where $\varrho_k$ and $\phi_k$, for all $k\in \mathcal E$, are assumed to define a steady-state solution consistent with specified values of compressor ratios  \cite{gugat2016stationary}.}  The network flow equations are defined by the initial-boundary-value system of PDEs \eqref{eq:net_flow}-\eqref{ic}.

%*********************************************************************************************************************************************************************
\section{Lumped Element Discretization} \label{sec:discretization}
%*********************************************************************************************************************************************************************

Spatial discretization of the network is formalized as a refinement of the graph topology. A refinement $(\hat{\mathcal E},\hat{\mathcal V})$ of the original graph $(\mathcal E, \mathcal V)$ is obtained by introducing auxiliary nodes, which are categorized as withdrawal nodes with zero outflow. These added nodes subdivide the edges of $\mathcal E$ so that the length of each refined segment satisfies $\ell_k\le \ell$ for all $k\in \hat{\mathcal{E}}$, where $\ell$ is a specified upper limit on the lengths of the pipeline segments in the refined graph.  Empirical studies have demonstrated that a segment length range of $5 \text{ km} \le \ell \le 10 \text{ km}$  is sufficient for both simulations and field data experiments in the physical regime of pipeline operations under usual conditions \cite{grundel2013computing}, \cite{zlotnik2017economic}. The refined graph inherits the prescribed edge orientations of the original graph.  For the sake of readability, hat symbols denoting the refined graph sets will be omitted in the following derivation.

Lumped element discretization is performed by integrating the dynamic equations \eqref{1}-\eqref{2} along the length of each pipeline segment $k\in \mathcal E$ so that \cite{himpe2021model} 
\begin{subequations} \label{eq:lumped_integration}
\begin{align}
\!\!\!\!\!\!\!\!\!\int_0^{\ell_k}\partial _t \rho_k +\partial_x \varphi_k dx & = 0,  \\
\!\!\!\!\!\!\!\int_0^{\ell_k}\delta\partial _t \varphi_k+ \sigma^2\partial_x\rho_k dx
& = -\frac{\lambda_k}{2D_k}\int_0^{\ell_k}\frac{ \varphi_k|  \varphi_k|}{\rho_k}dx. 
\end{align}
\end{subequations}
The above integrals of space derivatives are evaluated using the fundamental theorem of calculus.  Following a previous study \cite{himpe2021model}, the remaining integrals are evaluated by approximating density and mass flux over the relatively small refined segment with density at the outlet and mass flux at the inlet of the refined segment, respectively.  Factoring the approximate integrand terms out of the integrals and incorporating the actions of compressors in equations \eqref{eq:bc1}-\eqref{eq:bc2}, we obtain for each $(k:i\mapsto j)\in \mathcal E$ equations of the form
\begin{subequations} \label{eq:dis_net}
\begin{align}
\ell_k\dot {\bm \rho}_j &=\varphi_k^0-\varphi_k^{\ell}, \label{eq:dis1}
\\
\ell_k \delta \dot \varphi_k^0+\sigma^2 \left(\bm \rho_j-\mu_k\bm \rho_i \right) &= -\frac{\lambda_k \ell_k}{2D_k}\frac{ \varphi^0_k \left| \varphi^0_k \right|}{\bm \rho_j},  \label{eq:dis2}
\end{align}
\end{subequations}
where $\bm \rho_i=\bm s_i$ for $i\in \mathcal V_s$. A dot above a variable will always represent the time-derivative of the variable.

To obtain a matrix-vector formulation of the dynamical system, we introduce several matrices. Define the $E\times E$ diagonal matrices $L$, $K$, and $X$ with diagonal entries $L_{kk}=\ell_{k}$, $K_{kk}=\lambda_k/(2D_k)$, and $X_{kk}=\chi_k$.  Define the weighted edge-node incidence matrix $\Xi$ of size $E \times V$ component-wise by
\begin{equation}
\,\,\Xi_{ki}=  
\begin{cases}
-\mu_{k}(t), & \text{edge $k\in i_{\mapsto}$ leaves node $i$,}
\\
1, & \text{edge $k\in$$_{\mapsto}i$ enters node $i$},
 \\
0, & \text{else.}
\end{cases} \label{eq:xi}
\end{equation}
The rows of $\Xi$ correspond to the edges of the network and the columns correspond to the nodes. For each $(k:i\mapsto j)\in \mathcal E$, the $k$-th row of $\Xi$ contains exactly two nonzero entries that appear in the $i$-th and $j$-th columns with values $\Xi_{ki}=-\mu_k$ and $\Xi_{kj}=1$.   Define the $E\times V_s$ sub-matrix $N$ of $\Xi$ by the selection of columns $i\in \mathcal V_s$ and the $E\times V_w$ sub-matrix $M$ of $\Xi$ by the selection of columns $i\in \mathcal V_w$, where $V_s$ and $V_w$ denote the numbers of supply and withdrawal nodes, respectively.  The purpose of defining these sub-matrices is to separate the known nodal density variables at supply nodes from the unknown density variables at withdrawal nodes.  Define the $E\times V_w$ signed matrix $Q=\text{sign}(M)$ componentwise with the convention that $\text{sign}(0)=0$.  The signed matrix is well-defined by the lower bound on compressor variables.  Define the positive and negative parts of $Q$ by $Q_{\ell}$ and $Q_0$, respectively, so that $Q=Q_{\ell} +Q_0$ and $|Q|=Q_{\ell} - Q_0$, where $|A|$ denotes the componentwise absolute value of a matrix $A$.

Define inlet and outlet edge mass flux vectors by $ \varphi^0 =(\varphi_1^0,\dots, \varphi_E^0)'$ and $ \varphi^{\ell} =( \varphi_1^{\ell},\dots, \varphi_E^{\ell})'$.  Define the vector of supply node densities $\bm s=(\bm s_{1},\dots, \bm s_{V_s})'$, the vector of withdrawal node densities $\bm \rho=(\bm \rho_{V_s+1},\dots, \bm \rho_{V})'$, and the vector of withdrawal node outflows $\bm w =(\bm w_{V_s+1},\dots, \bm w_{V})'$, where the subscripts of the entries are indexed according to the node labels in $\mathcal V$.  
% Define also the vector of compressor ratios by $\mu = (\mu_1,\ldots,\mu_C)$, where $C=|\mathcal{C}|$.  
By applying the above matrix definitions, the discretized equations  \eqref{eq:dis_net} together with the boundary conditions \eqref{eq:boundary_conditions} can be written as a system of differential-algebraic equations (DAEs) of form
\begin{subequations} \label{eq:Dis_net}
\begin{align}
L Q_{\ell} \dot {\bm \rho} &= \varphi^0 -  \varphi^{\ell} , \label{eq:dDis1}
\\
\!\!\! L \delta \dot \varphi^0+\sigma^2 \left( M \bm \rho+N \bm s \right) &= -LK\frac{ \varphi^0 \odot | \varphi^0 |}{ Q_{\ell}\bm \rho}, \label{eq:dDis2}\\
\bm w&=Q_{\ell}' X\varphi^{\ell}+ Q_0'X\varphi^0, \label{eq:dDis3} 
\end{align}
\end{subequations}
where $\odot$ denotes component-wise multiplication, and the ratio of vectors on the right-hand side of equation \eqref{eq:dDis2} is defined component-wise. Henceforth, the ratio of two vectors will always denote component-wise division. 
\edit{Multiplying both sides of equation \eqref{eq:dDis1} by $ Q_{\ell}'X$, and then using \eqref{eq:dDis3} along with the relation $Q=Q_{\ell}+Q_0$, we obtain $Q_{\ell}' X L Q_{\ell}\dot {\bm \rho}=Q'X \varphi^0- \bm w$.}  Observe that the outlet flux becomes a dependent variable. Now, define the state vector of inlet flux by $\varphi=\varphi^0$.  Proposition 1 below provides a valuable result for inverting the so-called mass matrix $\Lambda=Q_{\ell}' X L Q_{\ell}$.   

{\bf Proposition 1} {\it The $V_w\times V_w$ matrix $\Lambda=Q_{\ell}' X L Q_{\ell}$ is diagonal with positive diagonal components
\begin{equation}
    \Lambda_{(j-V_s),(j-V_s)}=\sum_{k\in_{\mapsto}j}\chi_k\ell_k
\end{equation}
for each $j\in \mathcal V_w$.}

{\bf Proof}
% Let $j\in \mathcal V_w$ be an arbitrary withdrawal node. The columns of $Q_{\ell}$ and $(X L Q_{\ell})$ corresponding to $j\in \mathcal V_w$ (i.e. the $(j-V_s)$-th columns) each contain exactly $n_j$ nonzero entries with values $(Q_{\ell})_{kj}=1$ and $(X L Q_{\ell})_{kj}=\chi_k\ell_k$ that appear in the $k$-th entries corresponding to $k\in$$_\mapsto j$, where $n_j$ is the number of elements in $_\mapsto j$.  Since the graph of the network is connected, the number of edges entering node $j$ is strictly positive for each $j\in \mathcal V_w$.  All of the other elements in the $j$-th columns of the two matrices above are zero.  
First, suppose that $j\in \mathcal V_w$ is the only withdrawal node in the network.  Then $Q_{\ell}$ and $X L Q_{\ell}$ are vectors with nonzero entries in the $k$-th rows if and only if $k\in$$_{\mapsto}j$. Thus, $\Lambda=Q_{\ell}'(X L Q_{\ell})=\sum_{k\in_{\mapsto}j}\chi_k\ell_k$ is a scalar.  Now suppose $\mathcal V_w$ contains two or more withdrawal nodes and consider $i,j\in \mathcal V_w$.  If $i\not = j$, then the sets $_{\mapsto}i$ and $_{\mapsto}j$ are disjoint. Therefore, the column of $Q_{\ell}$ corresponding to $i\in \mathcal V_w$ and the column of $(X L Q_{\ell})$ corresponding to $j\in \mathcal V_w$ are orthogonal, so that $\Lambda_{(i-V_s),(j-V_s)}=\sum_{k\in \mathcal E}(Q_{\ell}')_{ik}(X L Q_{\ell})_{kj}=0$. If $i=j$, then again, the $(j-V_s)$-th columns of $Q_{\ell}$ and $X L Q_{\ell}$ have nonzero entries in the $k$-th rows if and only if $k\in$$_{\mapsto}j$, so
\begin{equation*}
    \Lambda_{(j-V_s),(j-V_s)}=\sum_{k\in \mathcal E}(Q_{\ell}')_{ik}(X L Q_{\ell})_{kj}=\sum_{k\in_{\mapsto}j}\chi_k\ell_k.
\end{equation*}
This completes the proof. \hfill $\square$

A practical consequence of Proposition 1 is that the inverse of $\Lambda$ can be computed by reciprocating its diagonal components. Geometrically, the $(j-V_s)$-th diagonal entry of $\Lambda$ corresponds to the total cylindrical volume of all refined pipeline segments directed to and incident upon the withdrawal node $j\in \mathcal V_w$. The shift by $V_s$ is related to the node indexing of $\mathcal V$ described above. Finally, the dynamical system of the pipeline network is given by
\begin{subequations} \label{eq:ode_system}
    \begin{align}
        \dot {\bm \rho}&=\Lambda^{-1}\left( Q'X \varphi- \bm w\right), \label{eq:ode_system1} \\
\delta \dot \varphi &= -\sigma^2L^{-1} \left( M \bm \rho+N \bm s \right)  -K\frac{\varphi \odot | \varphi |}{ Q_{\ell}\bm \rho}. \label{eq:ode_system2}
    \end{align}
\end{subequations}
The state vectors of the system include the withdrawal node density $\bm \rho$ and the edge inlet flux $\varphi$. The compressor variables are contained in the matrices $M$ and $N$. The system described by equation \eqref{eq:ode_system} is numerically consistent in the sense that, as the refined edge lengths approach zero, the discretized dynamics converge to the continuous dynamics defined by equations \eqref{eq:net_flow}–\eqref{eq:boundary_conditions} \cite{zlotnik2015optimal, baker2023optimal}. The initial condition for the discretized system corresponds to the sampled initial condition in equation \eqref{ic} at the nodes and edges of the refined network.  We define
\begin{equation}
 \bm \rho(0)=\bm \varrho, \qquad \varphi(0)=\phi, \label{eq:ic_ODE}
\end{equation}
and assume that $\varrho$ and $\phi$ represent a steady-state solution consistent with the initial values of compressor and boundary condition variables.

%*********************************************************************************************************************************************************************
\section{Adaptive Linear Control System}  \label{sec:linear_ti}
%*********************************************************************************************************************************************************************

In this section, we derive an adaptive linear system to approximate the nonlinear dynamics.  We suppose that such an approximation is made about a nominal state $\overline{\bm \rho}$ and $\overline \varphi$,  nominal boundary condition values $\overline{\bm s}$ and $\overline{\bm w}$, and nominal compressor-weighted matrices $\overline M$ and $\overline N$.  The nominal state is defined to be a currently estimated state and not necessarily a steady state of the system.  A linear model that governs the dynamics about an arbitrary nominal state provides flexibility to adaptively update the nominal state and linear system dynamics as time progresses to ensure the continued relevance of the linear model.  Linearization about constant nominal variables is motivated by a model-predictive control design applied on a moving horizon in which the nominal state, boundary condition values, compressor-weighted matrices, and corresponding linear dynamics are adaptively updated as time advances.

We linearize each term in equations \eqref{eq:ode_system1}–\eqref{eq:ode_system2} separately. Notably, equation \eqref{eq:ode_system1} is already linear with respect to the state, control, and boundary condition variables. Linearizing the friction term in equation \eqref{eq:ode_system2} about the nominal state with respect to $\bm \rho$ and $\varphi$ results in the $E\times E$ Jacobian sub-matrices
\begin{equation}\label{eq:alpha-beta-bar}
 \overline \alpha =K\text{diag}\left(\frac{ \overline \varphi \odot \left| \overline \varphi \right|}{ Q_{\ell} \overline{\bm \rho}^2 } \right), \quad 
 \overline \beta =2K\text{diag}\left(\frac{ \left| \overline \varphi \right|}{ Q_{\ell} \overline{\bm \rho }} \right), 
\end{equation}
respectively, where the square of $\overline{\bm \rho}$ is evaluated component-wise. As in Section \ref{sec:transfer_matrix}, the negative sign has been factored out of the $\overline \beta$ term. In comparison to the analysis above, the Jacobian sub-matrices $\overline \alpha$ and $\overline \beta$ are related to the coefficients $\alpha$ and $\beta$ in equation \eqref{eq:linear_dyn2}.  In particular, the set of coefficients $\alpha$ and $\beta$ in equation \eqref{eq:linear_dyn2} defined over each pipe segment of the refined network become the diagonal elements of the matrices $\overline \alpha$ and $\overline \beta$. By the chain rule, the total derivative of the friction term in equation \eqref{eq:ode_system2} with respect to $\bm \rho$ requires the post-multiplication of $\overline \alpha$ by $Q_{\ell}$. 

Linearizing equation \eqref{eq:ode_system} in $\bm \rho$, $\varphi$, $M$, $N$, $\bm s$, and $\bm w$ about the nominal state, control, and boundary condition variables results in the linear time-invariant system
\begin{equation}
\!\!\!\! \begin{bmatrix}
\dot{\bm \rho} \\
\delta \dot \varphi
\end{bmatrix}
 =\overline A \begin{bmatrix}
\bm \rho \\
\varphi
\end{bmatrix}
-
\begin{bmatrix}
\Lambda^{-1} \bm w \\
\sigma^2L^{-1}(N\overline{\bm s} + M\overline{\bm \rho}+ \overline N\bm s)
\end{bmatrix}
+\overline F,  \quad
\label{eq:linsys}
\end{equation}
where the state matrix $\overline A$ and additive nominal term $\overline F$ are given by
\begin{equation} \label{eq:state_matrix_A}
\overline A=
\begin{bmatrix}
0 & \Lambda^{-1}Q'X \\
\overline \alpha Q_{\ell}-\sigma^2 L^{-1} \overline M & -\overline \beta
\end{bmatrix}
\end{equation}
and
\begin{equation} \label{eq:F}
 \!\!\!\!  \overline F\!=\!\begin{bmatrix}
0 \\
 \sigma^2 L^{-1}\left( \overline M \overline{\bm \rho} +\overline N\overline{ \bm s} \right)
 \!-\! K\frac{\overline \varphi \odot \left| \overline \varphi  \right|}{ Q_{\ell} \overline{\bm \rho}} -\overline \alpha Q_{\ell}\overline{\bm \rho} +\overline \beta \overline \varphi
\end{bmatrix} \!.
\end{equation}
The initial condition of the linear system is equal to the initial condition in equation \eqref{eq:ic_ODE}.  If the nominal variables and linear system dynamics are updated as time progresses, then the initial condition would likewise be updated according to the current state. Note that the same symbols are used to describe both the nonlinear and linear system solutions, and we will continue to use these interchangeably when the context is clear. Before introducing the model-predictive controller, we first analyze the stability of the dynamics in equation \eqref{eq:ode_system} and the accuracy of the linear system approximation.

%*********************************************************************************************************************************************************************
\section{Eigenvalue and Pole Analysis}  \label{sec:eigenvalues}
%*********************************************************************************************************************************************************************

In this section, we examine the eigenvalues of $\overline A$ and the poles of the transfer matrix $G(\bm j \omega)$. Although pipeline network flows are stable from a control perspective because physical friction strongly dissipates kinetic energy, a stability analysis is still critical for designing controllers and preconditioners \cite{srinivasan2022numerical}.

%*********************************************************************
\subsection{Eigenvalues of the State Matrix}

We begin with the following proposition that establishes an identity on the eigenvalues of the state matrix $\overline A$.

{\bf Proposition 2}  {\it Let $\{\zeta(m)\}_{m=1}^{E+V_w}$ be the set of eigenvalues of $\overline A$ in equation \eqref{eq:state_matrix_A}, where $E$ and $V_w$ are the cardinalities of $\mathcal E$ and $\mathcal V_w$.  Then
\begin{equation}
    \sum_{m=1}^{E+V_w}\zeta(m) = -\sum_{k=1}^{E} \frac{\lambda_k |\overline \varphi_k|}{D_k(Q_{\ell}\overline{\bm \rho})_k}. \label{eq:center_mass}
\end{equation} }

{\bf Proof} The trace formula for the sum of eigenvalues of  $\overline A$ gives the result
\begin{equation*}
    \sum_{m=1}^{E+V_w}\zeta(m) = -\sum_{k=1}^{E} \overline \beta_{kk}=-\sum_{k=1}^{E} \frac{\lambda_k |\overline \varphi_k|}{D_k(Q_{\ell}\overline{\bm \rho})_k}. 
\end{equation*}
This concludes the proof. \hfill $\square$

We define the {\it center of gravity} of the eigenvalues of $\overline A$ as the average value $c=\sum_{m=1}^{E+V_w}\zeta(m)/(E+V_w)$.  Proposition 2 implies that the center of gravity is independent of $\overline \alpha$ and assumes a negative value for physical systems operating with positive densities.  The dependency of the center of gravity on $\overline \beta$ but not on $\overline \alpha$ aligns with the initial definition of the frequency domain representation in Section \ref{sec:transfer_matrix}. 
Interestingly, the center of gravity also does not explicitly depend on the network topology, described by the incidence matrices $\overline M$ and $\overline N$.  This suggests that the eigenvalues of the state matrix may be approximated by those associated with individual pipes. We revisit this hypothesis later in this section. First, we show that the state matrix associated with the simplified system, described by $\delta =0$ and $\overline \alpha=0$, is stable in the sense that the real parts of its eigenvalues are strictly negative.  Using the results from Section \ref{sec:pipe_flow}, let us suppose that the rates of variation are sufficiently small so that this simplified system may be considered.  

To simplify the algebra, we assume that the actions of compressors are unity in the following theorem.  With $\delta=0$ and $\overline \alpha=0$, we solve for $\varphi$ in the lower part of equation \eqref{eq:linsys} and substitute the resulting expression into the upper part of the linear system.  Because the weighted incidence matrix $\overline M$ equals the signed incidence matrix $Q$ when compressors are bypassed, the resulting state matrix of the simplified system is given by
\begin{equation} \label{eq:zero_input_response}
    \overline A=-\Lambda^{-1}(RQ)'RQ,
\end{equation}
where $R$ is the $E\times E$ diagonal matrix with positive diagonal components $R_{kk}=\sqrt{\sigma^2\pi D_k^3(Q_\ell\overline{\bm \rho})_k/(4\ell_k \lambda_k|\overline \varphi|)}$ for $k\in \mathcal E$.  We state a few preliminary results.   A real symmetric matrix has strictly positive eigenvalues if and only if it is positive definite \cite{khalil2002nonlinear}.  If $A$ and $B$ are arbitrary positive definite $n\times n$ matrices, then $ABA$ is positive definite \cite{horn2012matrix}. Suppose that $A$ and $B$ are positive definite, that $A$ is diagonal, and consider $AB=A^{1/2}(A^{1/2}BA^{1/2})A^{-1/2}$, where $A^{1/2}$ is the diagonal matrix with diagonal components equal to the square roots of those of $A$.  Because $A^{1/2}$ is positive definite, the above result says that $(A^{1/2}BA^{1/2})$ is positive definite.  Therefore, $AB$ is similar to a positive definite matrix and so its eigenvalues are real and strictly positive. 

{\bf Theorem 1} {\it Assume that the gas network is connected, has no self-loops, and its compressors are bypassed.  Then the $V_w\times V_w$ state matrix $\overline A=-\Lambda^{-1}(RQ)'RQ$ of the simplified linear system, described by $\delta=0$ and $\overline \alpha=0$, has strictly negative real eigenvalues.}

{\bf Proof}
As a consequence of Proposition 1, the diagonal matrix $\Lambda^{-1}$ is positive definite. A result from graph theory states that the number of linearly independent columns of $\Xi$ for a connected network with $V$ nodes is $(V-1)$ \cite{thulasiraman2011graphs}.  Because $Q$ is formed by removing $V_s\ge 1$ columns from $\Xi$, it follows that $Q$ has full column rank.  From Sylvester's rank inequality, the matrix $RQ$ also has full rank.  Thus, the matrix $(RQ)'RQ$ is invertible and, consequently, positive definite.  Applying the last of the results stated above shows that $\Lambda^{-1}(RQ)'RQ$ is similar to a positive definite matrix.  Therefore, we conclude that the eigenvalues of $\Lambda^{-1}(RQ)'RQ$ are real and positive, from which the result follows. 
\hfill $\square$

\begin{figure*}[!t]   
    \centering
    \begin{subfigure}[t]{0.4\textwidth}
        \centering
        \includegraphics[width=1\linewidth]{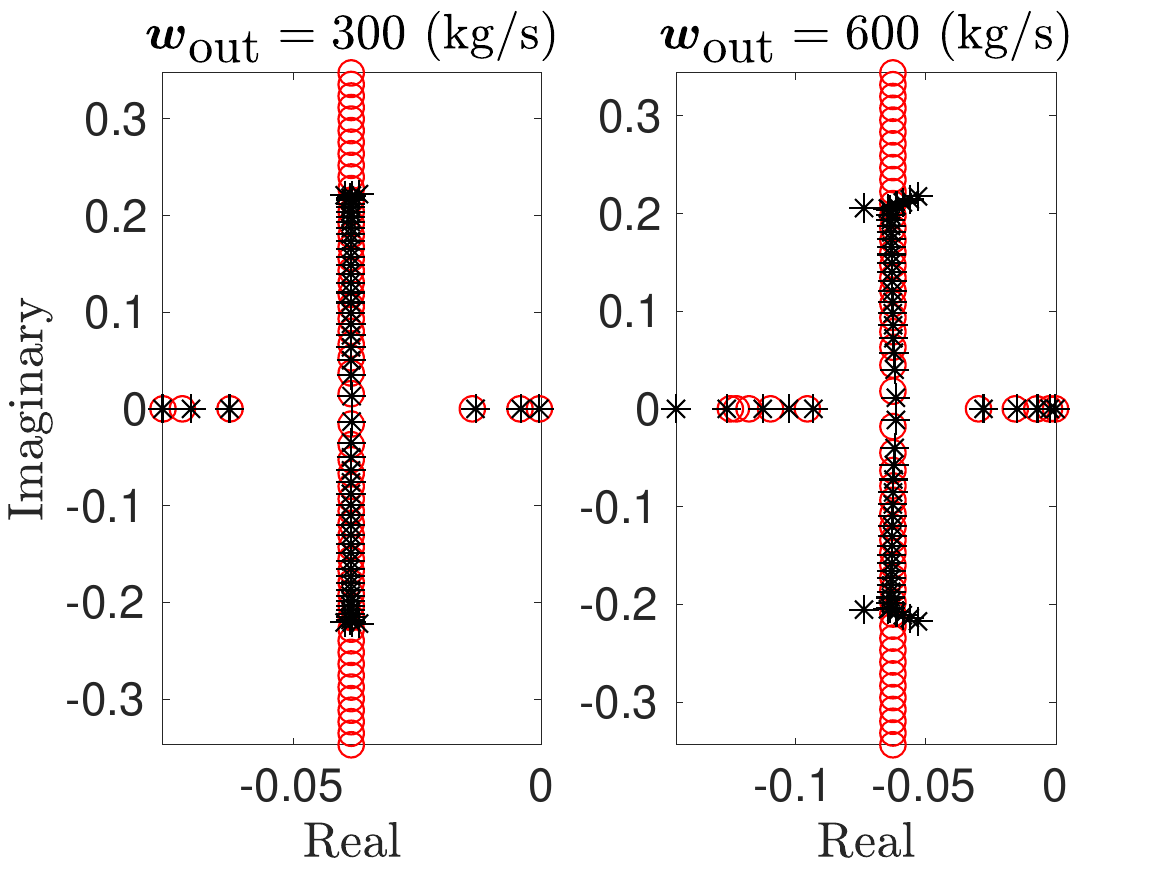}
        \caption{Single Pipe}
    \end{subfigure}%
    ~ 
    \begin{subfigure}[t]{0.4\textwidth}
        \centering
        \includegraphics[width=1\linewidth]{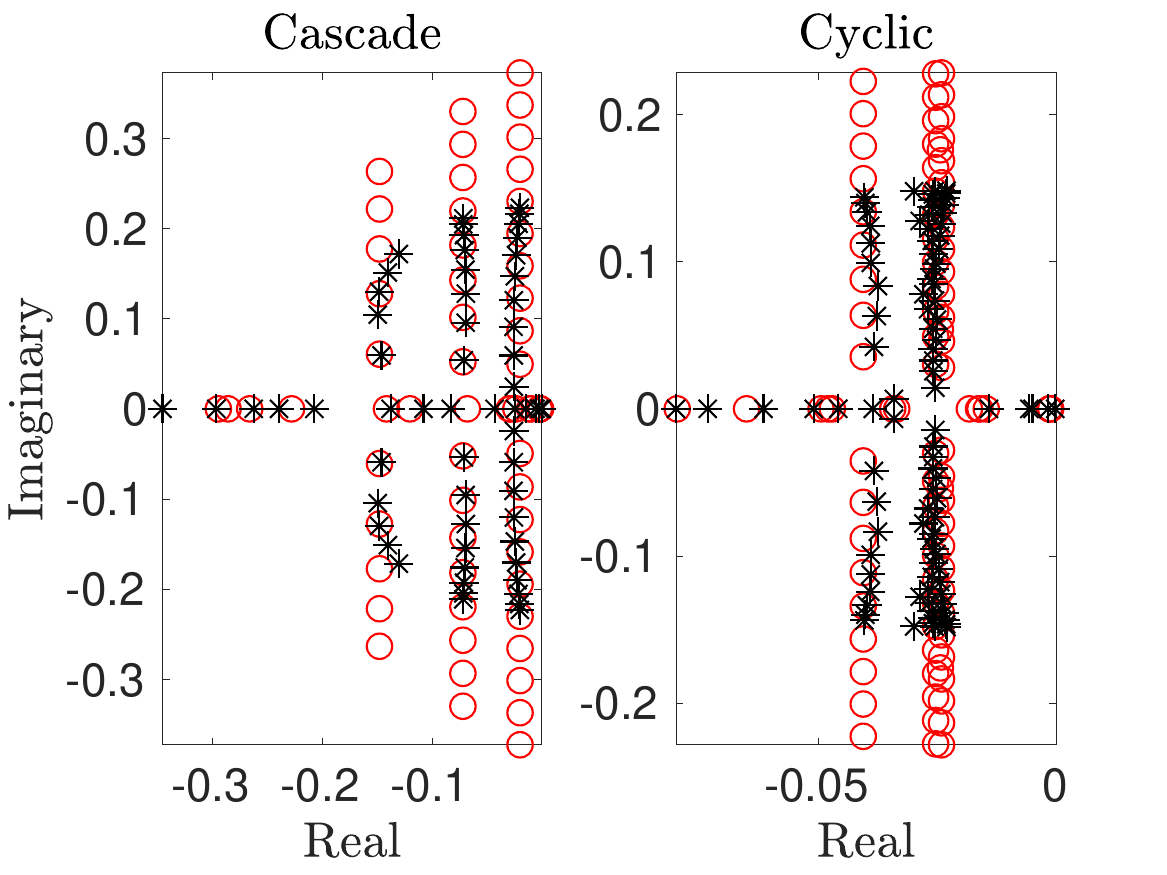}
        \caption{Networks}
    \end{subfigure}
    \caption{Eigenvalues ($\ast$) of the unsimplified state matrix $\overline A$ in \eqref{eq:state_matrix_A} and poles (\color{red}$\circ$\color{black})  of the transfer function $G(s)$ in \eqref{eq:transfer_matrix}. (a) A single pipe discretized into 30 segments and for two values of mass outflow (left and right). The pipeline parameters are $\ell=100$ km, $D=0.75$ m, $\lambda=0.01$, and $\sigma=377$ m/s. (b) The cascade network (left) consists of a series connection of three pipes with three withdrawal nodes before discretization. The cyclic network (right) is defined below in Figure \ref{fig:mpc_cyclic}.}
    \label{fig:eigA}
\end{figure*}

%*********************************************************************
\subsection{Poles of the Pipeline Transfer Matrix}

Although the eigenvalues of the finite-dimensional state matrix $\overline A$ for a single pipe in its refined edge representation generally differ from a subset of the infinitely many poles of the irrational transfer matrix $G(s)$, these poles have a well-defined closed-form expression, which can provide valuable insight into the properties of the eigenvalues of $\overline A$. Specifically, the poles of each component of $G(s)$ in equation \eqref{eq:transfer_matrix} for $\alpha = 0$ are given by 
\begin{equation}
    \zeta^{\pm}(m)=-\frac{\beta}{2}\pm \bm j \sqrt{\left(\frac{\pi \sigma}{2 \ell}\right)^2(2m+1)^2-\left(\frac{\beta}{2}\right)^2}, \label{eq:poles}
\end{equation}
for $m=0,1, \dots$.  Because the scalar $\beta$ in equation \eqref{eq:linear_dyn2} depends on the chosen nominal density, we consider the average value of the diagonal components of the matrix $\overline \beta$ and define $\beta=\sum_{k=1}^{E} \overline \beta_{kk}/E$.  From equation \eqref{eq:poles}, the imaginary asymptote of the poles of $G(s)$ intersects the real axis of the complex plane at $s=-\beta/2$.  The following proposition shows that this intersection point corresponds to the center of gravity.

{\bf Proposition 3} {\it Suppose that $\beta=\sum_{k=1}^{E} \overline \beta_{kk}/E$. Then the imaginary asymptote of the poles of $G(s)$ in equation \eqref{eq:transfer_matrix} for a single pipe with one supply node intersects the real axis of the complex plane at the center of gravity of the eigenvalues of $\overline A$ for the discretized pipe.}

{\bf Proof}
From Proposition 2, we have the relation $\beta=-c(E+V_w)/E$.  Because $E=V_w$ for a single discretized pipe with one supply node, it follows that $c=-\beta/2$. 
\hfill $\square$

Figure \ref{fig:eigA}(a) depicts the poles of $G(s)$ and the eigenvalues of $\overline A$ for a single pipe operating with two different withdrawal boundary conditions.  The effect of Proposition 3 is clearly evident. The figure further demonstrates that the poles and eigenvalues are qualitatively similar in that few are purely real while many are complex and collected along the imaginary asymptote.  Moreover, the poles of $G(s)$ provide quantitative approximations to the eigenvalues of $\overline A$, especially for small mass flux boundary condition values. In every case, all of the poles of $G(s)$ have strictly negative real parts for nonzero mass flux and positive density values, because $\beta$ in equation \eqref{eq:linear_dyn2} is always positive in such cases.  This observation agrees with the asymptotic stability of the system, as indicated by Theorem 1.

%*********************************************************************
\subsection{Poles of the Network Transfer Matrix}

Extending notation from Section \ref{sec:transfer_matrix}, let us define density and flow variations in the Laplace domain by $\bm P_k(s,x)$ and $\bm \Phi_k(s,x)$, respectively, for each edge $k\in \mathcal E$.  These variables are concatenated to form the vectors $\bm P=(P_1,\dots,P_E)'$ and $\bm \Phi=(\Phi_1,\dots,\Phi_E)'$.  Define the matrices $\bm \ell=\text{diag}(\ell_k)$, $\bm \gamma(s) =\text{diag}(\gamma_k(s))$, and $\bm z_c(s) =\text{diag}((z_c)_k(s))$.  Following a previous approach \cite{zecchin}, the single pipeline transfer matrix is extended to the network with the representation
\begin{equation} \label{eq:network_transfer_matrix}
\begin{bmatrix}
\bm P^{\ell} \\
\bm \Phi^{0}
\end{bmatrix}
=
\begin{bmatrix}
\text{sech} \left(\bm \ell \bm \gamma \right) & -\frac{\bm z_c}{\sigma} \text{tanh} \left(\bm \ell \bm \gamma  \right) \\
\frac{\sigma}{\bm z_c}\text{tanh}\left(\bm \ell \bm \gamma  \right) & \text{sech} \left(\bm \ell \bm \gamma \right)
\end{bmatrix}
\begin{bmatrix}
\bm P^{0} \\
\bm \Phi^{\ell}
\end{bmatrix},
\end{equation}
where we use the notation $g(\bm \ell \bm \gamma)=\text{diag}(g(\ell_k\gamma_k))$ for a complex-valued function $g$ of a complex argument.  For each $k\in \mathcal E$, the poles appearing in the $k$-th and $(E+k)$-th rows of the network transfer matrix are given by
\begin{equation}
    \zeta^{\pm}_k(m)\!=\!-\frac{\beta_k}{2}\pm \bm j \sqrt{\left(\frac{\pi \sigma}{2 \ell_k}\right)^2(2m+1)^2-\left(\frac{\beta_k}{2}\right)^2}, \label{eq:poles_network}
\end{equation}
for $m=0,1, \dots$.  This expression suggests that each pipe in the network may have its own imaginary asymptote in the complex plane, around which the complementing eigenvalues of $\overline A$ are clustered. Therefore, in a network with $E$ pipes, the eigenvalues may cluster around $E$ such asymptotes, which intersect the real axis of the complex plane at $c_k=-\beta_k/2$ for $k=1,\dots,E$. Clearly, the distance $|c_{k_1}-c_{k_2}|$ between any two asymptotes associated with two pipes is proportional to the absolute difference $|\beta_{k_1}-\beta_{k_2}|$ between their Jacobian terms.  A stronger interpretation arises when the pipes have the same diameters, friction factors, and operate at equal pressures. In this case, the definition of $\beta_k$ shows that the distance between the two asymptotes is proportional to the absolute difference between the associated mass flux values.  While diameters, friction factors, and pressures may vary across a network, we present this result because it provides some intuition for network operations.

To illustrate these effects, consider the same pipe in Figure \ref{fig:eigA}(a) with supply at the inlet and withdrawal at the outlet.  As in gunbarrel applications \cite{zhang2020transient}, we split this pipe into three adjacent pipes of equal lengths by adding two intermediate withdrawal nodes. From inlet to outlet, the three withdrawal nodes experience outflows of 500, 200, and 70 (kg/m$^2$s), respectively. From equation \eqref{eq:bc}, the corresponding mass flux values in these pipes are $\overline \varphi_1=770$, $\overline \varphi_2=270$, and $\overline \varphi_3=70$.  The eigenvalues of $\overline A$ and poles of the network transfer matrix are shown in the cascade plot on the left side of Figure \ref{fig:eigA}(b). Because the mass flux values differ significantly, the figure confirms the above result in that three distinct imaginary asymptotes appear in the complex plane. The right side of Figure \ref{fig:eigA}(b) displays the eigenvalues and poles corresponding to the 5-node cyclic network defined below in Figure \ref{fig:mpc_cyclic}.  This cyclic network contains five pipes, of which three deliver similar mass flux values. The above hypothesis is also evident in this example.

%*********************************************************************************************************************************************************************
\section{Error Analysis}  \label{sec:error}
%*********************************************************************************************************************************************************************

In this section, we examine the discrepancy between the approximate linearized dynamics \eqref{eq:linsys} and the nonlinear dynamics \eqref{eq:ode_system}. The main results are maximum bounds on the difference between the solutions of the linearized and nonlinear systems as functions of elapsed time and magnitude of variations about a nominal state. For simplicity, we consider practical pipeline settings and assume that $\bm s = \overline{\bm s}$, $M = \overline{M}$, and $N = \overline{N}$. Under these conditions, the linearized system in \eqref{eq:linsys} with $\delta = 1$ reduces to
 \begin{equation} 
 \!\!\!\begin{bmatrix}
\dot{\bm \rho} \\
\dot \varphi
\end{bmatrix}
\! =\! \overline A \begin{bmatrix}
\bm \rho \\
\varphi
\end{bmatrix}
\!+\!\begin{bmatrix}
-\Lambda^{-1} \bm w\\
 \!-\!\sigma^2 L^{-1}\overline N\overline{ \bm s}
 \!-\! K\frac{\overline \varphi \odot \left| \overline \varphi  \right|}{ Q_{\ell} \overline{\bm \rho}} \!-\!\overline \alpha Q_{\ell}\overline{\bm \rho}
 \!+\!\overline \beta \overline \varphi
\end{bmatrix}\!\!. \label{eq:state_matrix_A2}
\end{equation}
For the purpose of distinction, the solutions of the nonlinear system \eqref{eq:ode_system}-\eqref{eq:ic_ODE} and the linearized system \eqref{eq:ic_ODE}-\eqref{eq:linsys} with $\delta = 1$ are denoted by $(\bm \rho, \varphi)$ and $(\bm \rho_{\text{lin}}, \varphi_{\text{lin}})$, respectively. Both solutions are assumed to satisfy the same initial conditions as specified in equation \eqref{eq:ic_ODE}.  The error between the solutions is denoted by $\bm e=(\bm \rho,\varphi)-(\bm \rho_{\text{lin}},\varphi_{\text{lin}})$.  In the following, norms of vectors and matrices are defined with respect to the infinity norm.

{\bf Theorem 2}  {\it Assume that the real parts of the eigenvalues of $\overline A$ are negative; that $\bm s=\overline{\bm s}$, $M=\overline M$, and $N=\overline N$; that $\overline \varphi_k>0$ for all $k\in \mathcal E$; and that $\overline{\bm \rho}_j>0$ for all $j\in \mathcal V_w$.  If, for all $t\in [0,T)$, the solution of equation \eqref{eq:ode_system} satisfies $|\bm \rho_j(t) -\overline{\bm \rho}_j| \le \kappa \overline{\bm \rho}_j$ for all $j\in \mathcal V_w$ and $|\varphi_k(t)-\overline \varphi_k | \le \kappa \overline \varphi_k$ for all $k\in \mathcal E$, with $0\le \kappa \le \kappa_{\max}$ and $\kappa_{\max}<1$ sufficiently small, then there exist positive constants $a$ and $b$ such that } 
\begin{equation}
 \|\bm e(t) \|\le    \left( \frac{a\kappa^2}{(1-\kappa)^2}+ \frac{b\kappa^2}{1-\kappa}\right)\int_0^t\| e^{\overline A \tau} \|d\tau. \label{eq:bound}
\end{equation}
{\it The time horizon $T$ may be finite or infinite.}

{\bf Proof}  Because density and flux are assumed to be positive, the function $f(\bm \rho,\varphi,\bm w)$ on the right-hand side of equation \eqref{eq:ode_system} is continuously differentiable.  Invoking the mean value theorem, this function may be written as
 \begin{equation}  \label{eq:mean_value_dynamics}
\!\!\!f
\! =\!A^* \! \begin{bmatrix}
\bm \rho \\
\varphi
\end{bmatrix}
\!+\!\begin{bmatrix}
-\Lambda^{-1} \bm w\\
\! -\sigma^2 L^{-1}\overline N\overline{ \bm s}
 \!-\! K\frac{\overline \varphi \odot \left| \overline \varphi  \right|}{ Q_{\ell} \overline{\bm \rho}} \! -\!\alpha^* Q_{\ell}\overline{\bm \rho}
 \!+\!\beta^* \overline \varphi
\end{bmatrix}\!\!, \! 
\end{equation}
where we define the matrices $\alpha^*=K\text{diag}(\varphi^*/ Q_{\ell} \bm \rho^*)^2$, $\beta^* =2K\text{diag}( \varphi^* / Q_{\ell} \bm \rho^*)$, and
\begin{align}
\!\!\! A^*&\!=\!
\begin{bmatrix}
0 & \Lambda^{-1}Q'X \\
\alpha^* Q_{\ell}-\sigma^2 L^{-1} \overline M & -\beta^*
\end{bmatrix}.  \label{eq:state_matrix_A_star}
\end{align}
Each component of the vector $(\bm \rho^*,\varphi^*)$ lies on the line segment joining the corresponding components of the vectors $(\bm \rho,\varphi)$ and $(\overline{\bm \rho},\overline \varphi)$. Using equations \eqref{eq:state_matrix_A2} and \eqref{eq:mean_value_dynamics}, the error dynamics are governed by
\begin{align}
    \dot{\bm e} 
=\overline A \bm e+\left( A^*-\overline A \right)
\begin{bmatrix}
\bm \rho -\overline{\bm \rho}\\
\varphi-\overline \varphi
\end{bmatrix}, \label{eq:error_dynamics}
\end{align}
where we use the fact that 
\begin{equation}
    A^* - \overline A =
\begin{bmatrix}
0 & 0 \\
(\alpha^*-\overline \alpha) Q_{\ell} & -(\beta^*-\overline \beta)
\end{bmatrix}. \label{eq:state_matrix_difference}
\end{equation}
% Multiplying by $P$ and applying a change of variables, the relative error satisfies 
% \begin{equation} \label{eq:rel_error_dynamics}
%     \dot{\bm e}_{\text{rel}}  =P\overline AP^{-1} \bm e_{\text{rel}}+P\left( A^*-\overline A \right)
% \begin{bmatrix}
% \bm \rho -\overline{\bm \rho}\\
% \varphi-\overline \varphi
% \end{bmatrix}.
% \end{equation}
\edit{The fundamental solution of equation \eqref{eq:error_dynamics} with zero initial error yields the general bound
\begin{equation} \label{eq:tight_error_bound}
  \!\!\!  \| \bm e(t) \| \le \int_0^t  \left\Vert e^{\overline A \tau} (A^*-\overline A)  \begin{bmatrix}
\bm \rho -\overline{\bm \rho}\\
\varphi-\overline \varphi
\end{bmatrix} \right\Vert  d\tau. 
\end{equation}
The sub-multiplication property of the norm will be applied.}  Using \eqref{eq:state_matrix_difference}, the difference between state matrices simplifies to differences between Jacobian sub-matrices. The $k$-th diagonal entries of these differences, with $k:i\mapsto j$, are given by
\begin{align*}
\alpha_{kk}^*-\overline \alpha_{kk} &= K_{kk}\left(\frac{ (\varphi_k^*)^2 }{(\bm \rho_{j}^* )^2} 
-\frac{ ( \overline \varphi_k )^2}{\overline {\bm \rho}_{j}^2}\right), \\
\beta_{kk}^*-\overline \beta_{kk} &= 2K_{kk}\left(\frac{\varphi_k^*}{\bm \rho_{j}^* } 
-
\frac{ \overline \varphi_{k} }{\overline{ \bm \rho}_{j}} \right).
\end{align*}
Because $\|\bm\rho-\overline{\bm \rho}\| \le \kappa \|\overline{\bm\rho}\|$ and $\|\varphi-\overline \varphi \| \le \kappa \|\overline \varphi\|$, it can be shown that
\begin{equation*}
\| \alpha^*-\overline \alpha \| \le  \frac{4\kappa}{(1-\kappa)^2}\|\overline \alpha \|, \quad 
\| \beta^*-\overline \beta \| \le \frac{2\kappa}{1-\kappa} \| \overline \beta \|.
\end{equation*}
Therefore,
\begin{equation*}
\!\!\! \left\Vert (A^*-\overline A)\begin{bmatrix}
\bm \rho -\overline{\bm \rho}\\
\varphi-\overline \varphi
\end{bmatrix} \right\Vert 
  \le  \frac{a \kappa^2}{(1-\kappa)^2}+ \frac{b\kappa^2}{1-\kappa},
\end{equation*}
where $a=4 \| \overline \alpha \|  \| \overline{\bm \rho} \|$ and $b=2 \| \overline \beta \| \| \overline \varphi \|$. The result follows by applying this inequality in \eqref{eq:tight_error_bound}. \hfill $\square$

\begin{figure}
\centering
\includegraphics[width=1.0\linewidth]{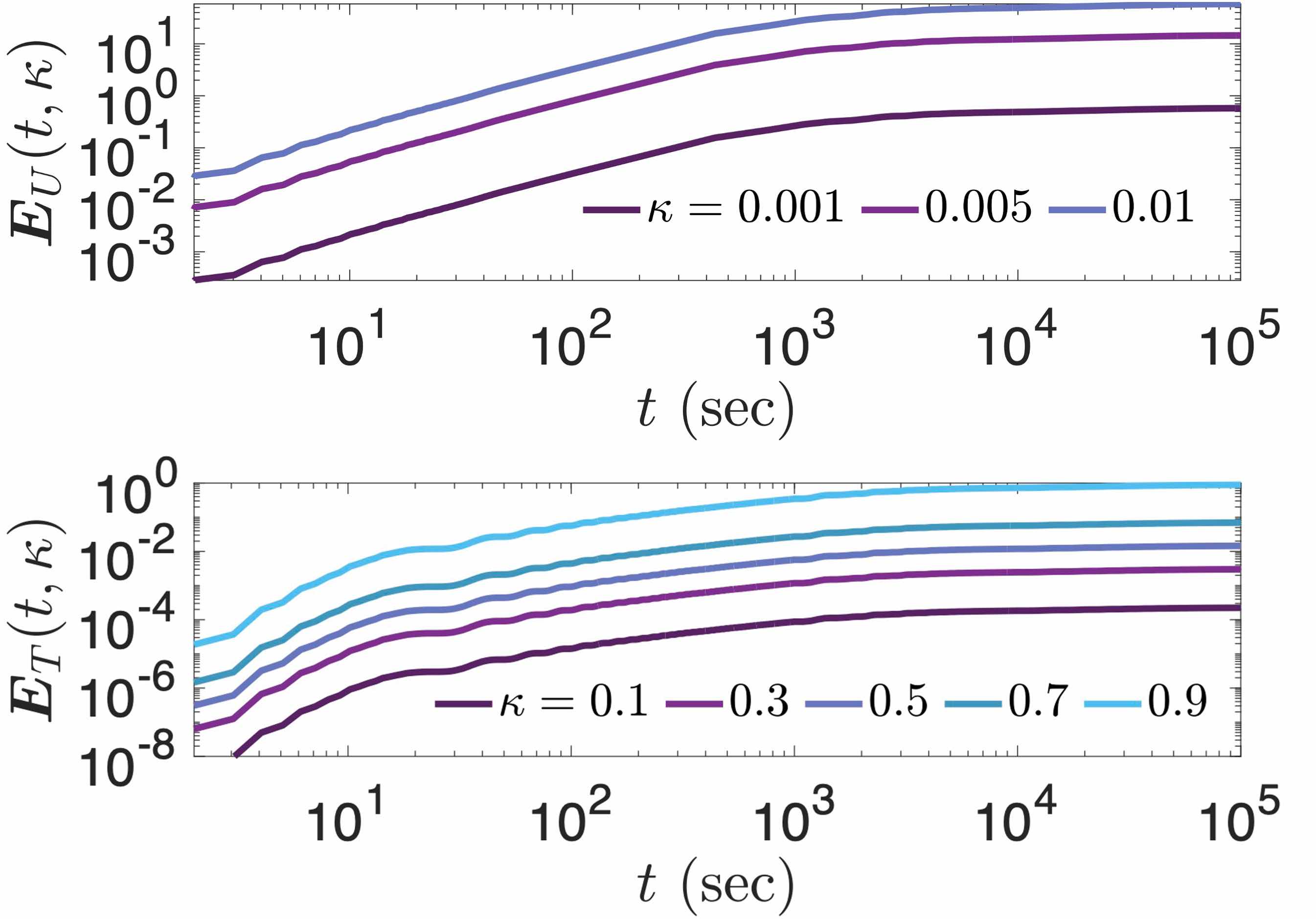}
\caption{\edit{Error bounds in equations \eqref{eq:uniform_error_bound} and \eqref{eq:tv_error_bound} for various values of $\kappa$ using the 25-node network example in Section \ref{sec:results}.}}
\label{fig:error_bound}
\end{figure}

In the proof of Theorem 2, we assume that $\kappa\le \kappa_{\text{max}}$ to ensure that the nonlinear solution is well-defined. Calibration of maximal appropriate values of $\kappa_{\text{max}}$  depends on empirical observations and is beyond our scope.  \edit{Using \eqref{eq:bound}, we define the \emph{uniform error bound} relative to the nominal state by
\begin{equation}
 \bm E_U(t,\kappa) =  \left( \frac{a_0\kappa^2}{(1-\kappa)^2}+ \frac{b_0\kappa^2}{1-\kappa}\right)\int_0^t\|e^{\overline AP \tau} \|d\tau,  \label{eq:uniform_error_bound}
\end{equation}
where $a_0=a/\|(\overline{\bm \rho},\overline \varphi)\|$ and $b_0=b/\|(\overline{\bm \rho},\overline \varphi)\|$.  
Observe that $\bm E_U(0,\kappa)=0$ and, for stable state matrices, $\bm E_U(t,\kappa)$ monotonically increases and converges as $t\rightarrow \infty$ for each value of $\kappa \le \kappa_{\max}$. This monotone convergence is illustrated in Figure \ref{fig:error_bound} for several values of $\kappa$. The state matrix used for this demonstration is obtained from the 25-node network example defined below in Section \ref{sec:results}. Although the uniform bound may have some practical implications, it does not account for the rates of variation nor the specific topological locations where these variations occur. If additional information is available, then the above proof may be modified to potentially provide a less conservative bound.

{\bf Theorem 3}  {\it Assume the conditions in Theorem 2 hold. Let $0\le \gamma(t)\le \gamma_{\max}$ for all $t\in[0,T)$, where $\gamma_{\max}<1$.  If, for all $t\in [0,T)$, the solution of equation \eqref{eq:ode_system} satisfies $|\bm \rho_j(t) -\overline{\bm \rho}_j| \le \gamma(t) \overline{\bm \rho}_j$ for all $j\in \mathcal V_w$ and $|\varphi_k(t)-\overline \varphi_k | \le \gamma(t) \overline \varphi_k$ for all $k\in \mathcal E$, then there exists a positive constant $c$ such that}
\begin{align}
\|\bm e(t) \|&\le    c\int_0^t \gamma^2(\tau) \| e^{\overline A \tau} \|d\tau. \label{eq:tv_bound}
\end{align}

{\bf Proof} Following the proof of Theorem 2, we have
\begin{equation*}
\!\!\! \left\Vert (A^*-\overline A)\begin{bmatrix}
\bm \rho -\overline{\bm \rho}\\
\varphi-\overline \varphi
\end{bmatrix} \right\Vert
   \le  \frac{a\gamma^2}{(1- \gamma_{\max})^2}  + \frac{b\gamma^2}{1-\gamma_{\max}} ,
\end{equation*}
where $a$ and $b$ are defined above.  The result follows from \eqref{eq:tight_error_bound} with $c=a/(1-\gamma_{\max})^2+b/(1-\gamma_{\max})$. 
% Using this inequality, the relative error is bounded by
% \begin{align*} 
%  \!\!\! \| \bm e_{\text{rel}}(t) \| &\le \left(\frac{a_0\kappa^2}{(1-\kappa)^2}  + \frac{b_0\kappa^2}{1-\kappa} \right)\int_0^t  \| e^{P\overline A P^{-1} \tau}\| \gamma^2(\tau)d\tau. 
% \end{align*}
% The desired bound is obtained by applying H{\"o}lder's inequality, where $a=a_0/\sqrt{4p+1}$ and $b=b_0/\sqrt{4p+1}$.
\hfill $\square$

The bound on the right-hand side of \eqref{eq:tv_bound} is called the \emph{time-varying error bound}. As with the uniform bound, the general result in Theorem 3 is illustrated for the 25-node network example defined below.  We now suppose that this network is operating with flow variations bounded by sinusoids of 1 cyc/hr, characterized by $\gamma(t)=\kappa|\sin(2\pi t/T)|$, where $T=3600$ sec. By substituting this expression into \eqref{eq:tv_bound} and normalizing the error, we obtain the relative error bound defined by
\begin{equation}
    \bm E_{T}(t,\kappa)=c_0\int_0^t \kappa^2 \sin^2(2\pi \tau/T)  \| e^{\overline A \tau} \|d\tau, \label{eq:tv_error_bound}
\end{equation}
where $c_0=c/\|(\overline{\bm \rho},\overline \varphi)\|$.  Figure \ref{fig:error_bound} displays this bound as a function of time for various values of $\kappa$.  Although the values of $\kappa$ are substantially large, the corresponding error bounds are mostly moderate and generally reflect the types of errors that are commonly observed in simulations \cite{baker2021analysis}.  

As a practical application, the time-varying error bound may be inverted numerically to determine necessary frequencies of a feedback controller that ensure the linear system dynamics in \eqref{eq:linsys} adequately resemble the nonlinear dynamics in \eqref{eq:ode_system}.  To demonstrate this process, suppose that the control operators of the 25-node network tolerate relative errors of at most 1\%.  If flow variations are bounded according to $\gamma(t)$ defined above with $\kappa=0.3$, then Figure \ref{fig:error_bound} suggests that neither the nominal state nor the linear model would require updates to maintain the permitted error tolerance.  However, if $\kappa=0.7$, then the operators should expect to update the nominal state, initial condition, and linear system model approximately every half hour in order to ensure that the linear model is accurate to within 1\% relative error.  This example inspires a model-predictive control design for compressor management in which initial conditions and dynamics are adaptively updated as the time window advances.}

%*********************************************************************************************************************************************************************
\section{Model-Predictive and Optimal Control} \label{sec:MPC}
%*********************************************************************************************************************************************************************

Building on previous studies \cite{abbaspour2008nonisothermal, zlotnik2015optimal}, we design a model-predictive control synthesis that seeks to minimize compressor energy consumption while ensuring that pressure, compression, and flow constraints are satisfied.  The time interval $[0,T]$ is discretized into small segments $(t_{m-1},t_m]$, for $m=1,\dots,m_T$, where each $t_m$ represents a sampling time at which measurements can be made.  The MPC optimizes flows on $(t_{m-1},t_m]$, and uses the resulting predictions to inform optimization over the next segment $(t_{m},t_{m+1}]$. This process is repeated recursively until control actions are determined for the entire time horizon $[0,T]$. To demonstrate the effectiveness of adaptive linear systems, we introduce two MPC designs.  One uses the nonlinear system \eqref{eq:ode_system} as its dynamic constraints, and the other uses the adaptive linear approximation \eqref{eq:linsys} that updates the nominal state and linear system matrices at each sampling time.  Additionally, we also synthesize an optimal control design for comparison to an optimal baseline solution.  This controller performs a single dynamic optimization across all time samples \cite{zlotnik2015optimal}.  \edit{When MPC and OC designs are implemented with nonlinear objectives and nonlinear dynamic constraints, we refer to them as nonlinear MPC and nonlinear OC, respectively. Similarly, when the objective functions and constraints are linear, we refer to them as linear MPC and linear OC.}

%*********************************************************************
\subsection{Model-Predictive Control}

For the current time segment $(t_{m-1},t_{m}]$, the nonlinear MPC performs an internal optimization to minimize the cost objective function defined by \cite{wong1968optimization, zlotnik2015optimal}
\begin{equation}
\mathcal J(t_m) \!=\!\sum_{k\in \mathcal C}c_k \varphi_{k}(t_m)\left(  \mu_{k}(t_m)^{(\gamma-1)/\gamma}-1 \right), \label{J}
\end{equation}
where $c_{k}$ is related to the efficiency of the compressor $k\in\mathcal{C}$, and $\gamma$ is the isentropic exponent of natural gas \cite{ivan2005}. 
For the linear MPC, we define the linearized objective function by 
\begin{align} \label{Jlin}
\mathcal J(t_m) =&\sum_{k\in \mathcal C}c_k \varphi_k(t_m) \left(  \overline \mu_k^{(\gamma-1)/ \gamma}-1\right) \nonumber \\
&+ \sum_{k\in \mathcal C}c_k
\left(\frac{\gamma-1}{\gamma } \overline \varphi_k\overline \mu_k^{-1/\gamma}  \right) \mu_k(t_m).
\end{align}
Pressure, flow, and compressor limitations are defined by inequality constraints of the form
\edit{\begin{equation}
\!\!\! \begin{bmatrix}
    \bm \rho^{\min} \\
    \varphi^{\min}
\end{bmatrix}   
\le   
\begin{bmatrix}
    \bm \rho (t_m) \\
    \varphi(t_m)
\end{bmatrix}  \le 
 \begin{bmatrix}
    \bm \rho^{\max} \\
    \varphi^{\max}
\end{bmatrix} , \; 1 \le   \mu_{k}(t_{m}) \le \mu_k^{\max}, \label{eq:ineq}
\end{equation}
where bounding inequalities are applied componentwise.} \revise{While nonlinear operation envelopes offer a more realistic representation of compressor unit limitations \cite{wu2000model}, we adopt box constraints in this study to remain consistent with prior simulation studies \cite{wong1968optimization, de2000gas, borraz2013optimization, misra2014optimal, zlotnik2015optimal}. Nonetheless, let us note that each of the proposed methods can accommodate more accurate representations by replacing box constraints on compressor variables with linear inequality constraints that represent hyperplane approximations of the envelopes \cite{wu2000model}.}  By implicitly integrating the continuous dynamics over the time segment $(t_{m-1},t_{m}]$, we further constrain the solution of the internal optimization to satisfy the discretized dynamics defined by
\begin{equation}
    \begin{bmatrix}
    \bm \rho (t_m) \\
    \varphi(t_m)
\end{bmatrix}  
=\begin{bmatrix}
    \overline{\bm \rho}  \\
    \overline \varphi
\end{bmatrix}+\Delta t_m \overline f(\bm \rho(t_m),\varphi(t_m),\mu_k(t_m)), \label{eq:OCP_dynamics}
\end{equation}
where $\Delta t_m=t_{m}-t_{m-1}$ and $\overline f$ represents the right-hand side of either the nonlinear system \eqref{eq:ode_system} or the updated linear system \eqref{eq:linsys}.  

The solution $\bm \rho^*(t_m)$, $\varphi^*(t_m)$, and $\mu_k^*(t_m)$, for $k\in \mathcal C$, is determined by solving the optimization problem
\begin{subequations} \label{eq:feed_ocp}
  \begin{alignat}{3}
&\text{minimize} & \qquad &  \text{Compressor energy: \eqref{J} or \eqref{Jlin},}\\
&\text{subject to}   & \qquad &\text{Inequality constraint: \eqref{eq:ineq},}\\
&& \qquad & \text{Dynamic constraint: \eqref{eq:OCP_dynamics}.} 
\end{alignat}
\end{subequations}
The solution is used to update the nominal variables according to $\overline{\bm \rho}=\bm \rho^*(t_m)$, $\overline \varphi=\varphi^*(t_m)$, and $\overline \mu_k=\mu_k^*(t_m)$. The process is then repeated for the next time segment $(t_m,t_{m+1}]$.  The linear and nonlinear MPCs are initialized by setting the nominal state as the initial condition defined in \eqref{eq:ic_ODE}.  \edit{Density, mass flux, and compressor ratios are all treated as decision variables in the above optimization problems. The inequality and dynamic constraints of the optimization problems are considered feasible if there exist values for the decision variables that simultaneously satisfy all of these constraints. Provided that the constraints remain feasible, one of the key advantages of linear MPC is that the convexity of its internal optimization ensures that any local optimal solution is globally optimal \cite{mayne2000constrained}.}

%*********************************************************************
\subsection{Optimal Control}

The objective function for OC is defined as \cite{zlotnik2015optimal}
\begin{equation} \label{eq:cumulative_energy}
   J = \sum_{m=1}^{m_T} \mathcal J(t_m),
\end{equation}
where $\mathcal J(t_m)$ is given by either equation \eqref{J} or \eqref{Jlin}, depending on whether linear or nonlinear OCs are formulated.  As with MPC, the linear and nonlinear OCs are both implemented in order to compare the outcomes of linear and nonlinear representations as the dynamic constraints.  The dynamic constraints of the linear and nonlinear OCs are defined to be the collection of equations \eqref{eq:OCP_dynamics}, for all $m=1,\dots,m_T$, derived using either a fixed linear system \eqref{eq:linsys} or the nonlinear system \eqref{eq:ode_system}, respectively.

The OC design optimizes policies by accounting for variations in boundary conditions that are anticipated to occur over the entire time horizon $[0,T]$.  In contrast, the proposed MPC only considers predicted boundary conditions for the next immediate time segment.  Consequently, we expect the MPC solution to be sub-optimal with respect to the OC solution for both the nonlinear and linear programs.  However, there are two potential advantages to using MPC.  First, it can readily adapt to real-time deviations from anticipated loads caused by disturbances or intra-day contract adjustments.  Second, the total computational cost of MPC may be orders of magnitude lower than that of OC.  Indeed, a worst-case computational complexity for a linear optimization program has been shown to scale as $\mathcal O(n^{2.5})$, where $n$ is the number of optimization variables \cite{vaidya1989speeding}. \edit{We use this result to establish a computational reduction provided by MPC.

{\bf Theorem 4}  {\it In the worst case scenario, the computational cost of linear OC scales $m_T^{1.5}$ times larger than that of linear MPC, where $m_T$ is the number of time samples.}

{\bf Proof} The number of optimization variables for OC is $n=m_T(E+V_w+C)$, where $C$ is the number of control variables.  Based on the result referenced above, the computational cost associated with linear OC scales as $\mathcal O(m_T^{2.5}(E+V_w+C)^{2.5})$.  On the other hand, the number of optimization variables for MPC defined on a single time segment is $n=(E+V_w+C)$.  Hence, the computational cost for linear MPC over a single time segment is of order $\mathcal O((E+V_w+C)^{2.5})$.  Accumulating costs for each segment gives the total computational cost that scales as $\mathcal{O}(m_T(E+V_w+C)^{2.5})$. It is worth noting that $m_T^{1.5}$ already exceeds two orders of magnitude for just hourly samples over a 24-hour time interval.
\hfill $\square$}

%\edit{We note that because the objective function in equation \eqref{J} is non-convex and the nonlinear dynamic constraints \eqref{eq:ode_system} are nonlinear, the nonlinear MPC and OC formulations may have multiple local minima.  A major advantage of the linearized MPC formulation that uses the objective in equation \eqref{Jlin} and the linearized dynamics in equation \eqref{eq:linsys} is convexity of the objective function.  A comprehensive analysis to characterize the uniqueness of solutions to these formulations is beyond the scope of our study.}

%*********************************************************************************************************************************************************************
\section{Computational Examples}  \label{sec:results}
%*********************************************************************************************************************************************************************

This section presents several numerical examples using two network topologies to demonstrate the effectiveness of linear and nonlinear MPC in comparison to linear and nonlinear OC.  Specifically, we assess the performance of these four control designs based on computation time and compressor energy accumulated over the entire operational horizon, as defined in \eqref{eq:cumulative_energy}.  To evaluate the applicability of linear MPC and linear OC, we measure the solution errors of these programs relative to their nonlinear counterparts. Rather than using the maximum relative error, as defined in Section \ref{sec:error}, we examine relative errors in density, flux, and control actions individually. For example, the maximum relative error in density is defined as
\begin{equation}
    \bm E_{\bm \rho}=2\max_m \frac{\|\bm \rho_{\text{non}}(t_m)-\bm \rho_{\text{lin}}(t_m)\|}{\| \bm \rho_{\text{non}}(t_m)+\bm \rho_{\text{lin}}(t_m)\|}\times 100\%,
\end{equation}
where, as above, we use maximum norms.  Similar definitions are made for flux and control errors, and are denoted by $\bm E_{\varphi}$ and $\bm E_{\mu}$, respectively.  For each of the two network topologies, we consider a standard 24-hr operational window $[0,T]$, where $T=24 \times 3600$s. In this analysis, several sets of equally spaced sampling times are applied for both MPC and OC. The uniform sampling interval $\Delta t$ is emphasized by referring to the controllers as MPC($\Delta t$) and OC($\Delta t$). All computations are performed in Matlab R2023a on a MacBook Pro. 

\subsection{5-Node Cyclic Network}

\begin{figure*}[t]
    \centering
    \begin{subfigure}[t]{0.333\textwidth}
        \centering
        \includegraphics[width=1\linewidth]{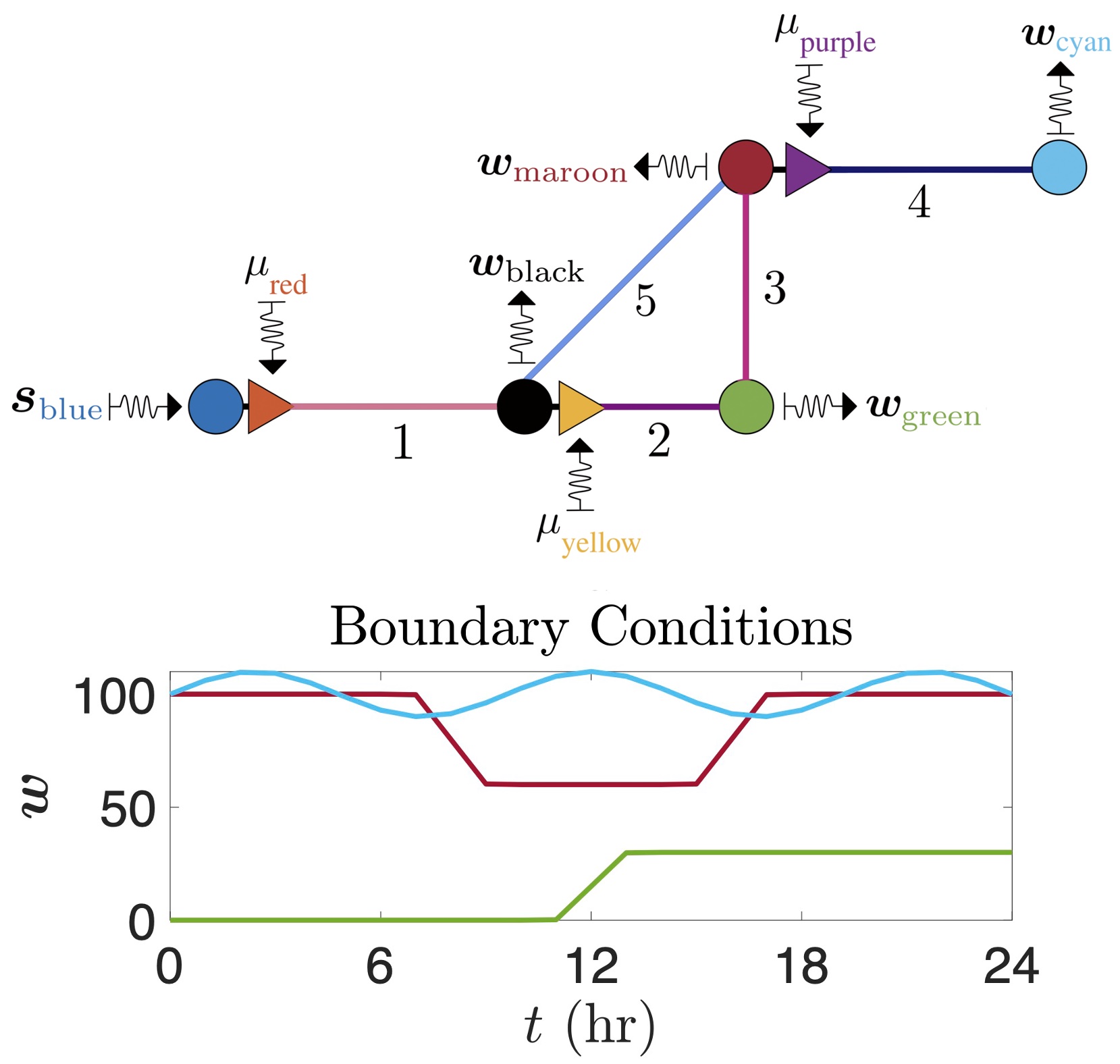}
    \end{subfigure}%
    ~ 
    \begin{subfigure}[t]{0.333\textwidth}
        \centering
        \includegraphics[width=1\linewidth]{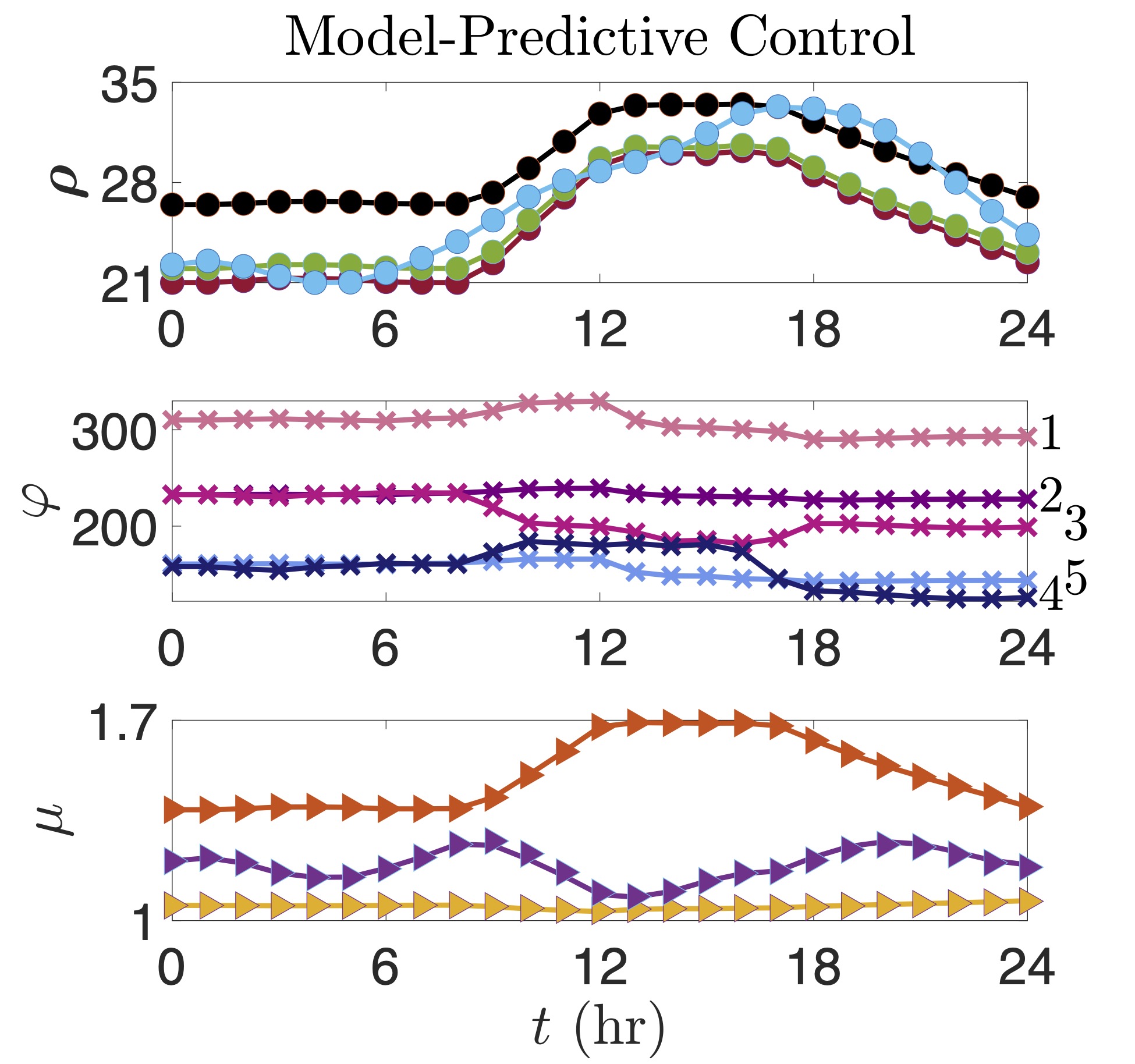}
    \end{subfigure}%
    ~ 
    \begin{subfigure}[t]{0.333\textwidth}
        \centering
        \includegraphics[width=1\linewidth]{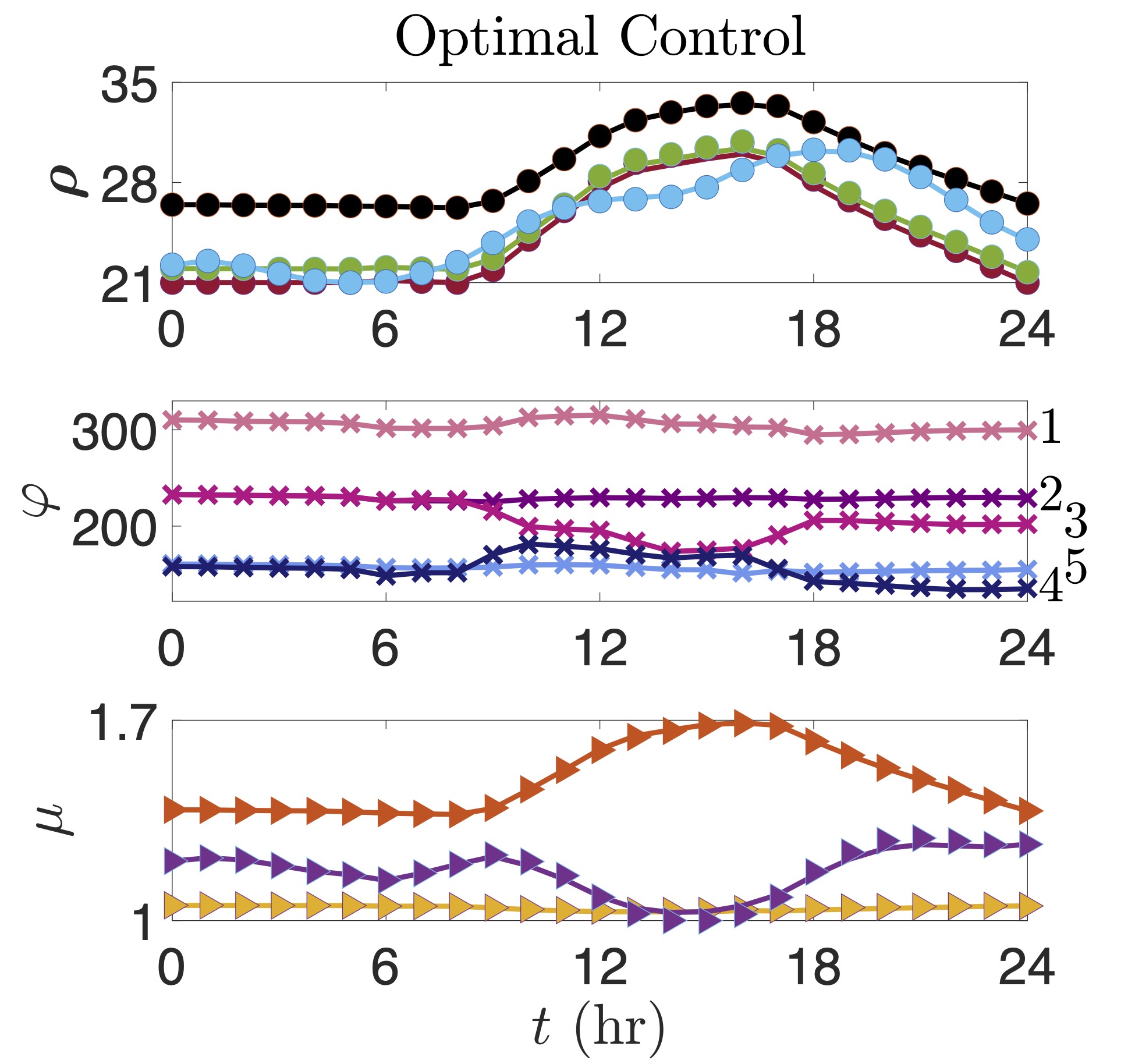}
    \end{subfigure}
    \caption{Left: Network configuration (not to scale) and boundary conditions at color-coordinated nodes \cite{grail2018data}.  Lines, circles, and triangles in the network diagram represent pipelines, nodes, and compressors, respectively. Pipe lengths in km are $\ell_1=20$, $\ell_2=70$, $\ell_3=10$, $\ell_4=80$, $\ell_5=60$.  The pipes have diameter $D_k=0.9144$ m and friction factor $\lambda_k=0.01$ for $k=1$, 2, 3, and 4, except for edge $5\in \mathcal E$ for which $D_5=0.635$ m and $\lambda_5=0.015$. The sound speed is $\sigma=377.964$ m/s. Center: MPC(60 min) solution. Right: OC(60 min) solution. Solid lines and markers represent solutions obtained using nonlinear and linear dynamics, respectively. The colors of density, flux, and compressor ratio plots correspond to respective locations in the network.}
    \label{fig:mpc_cyclic}
\end{figure*}

\begin{figure*}[ht]
    \centering
    \begin{subfigure}[t]{0.333\textwidth}
        \centering
       \includegraphics[width=1\linewidth]{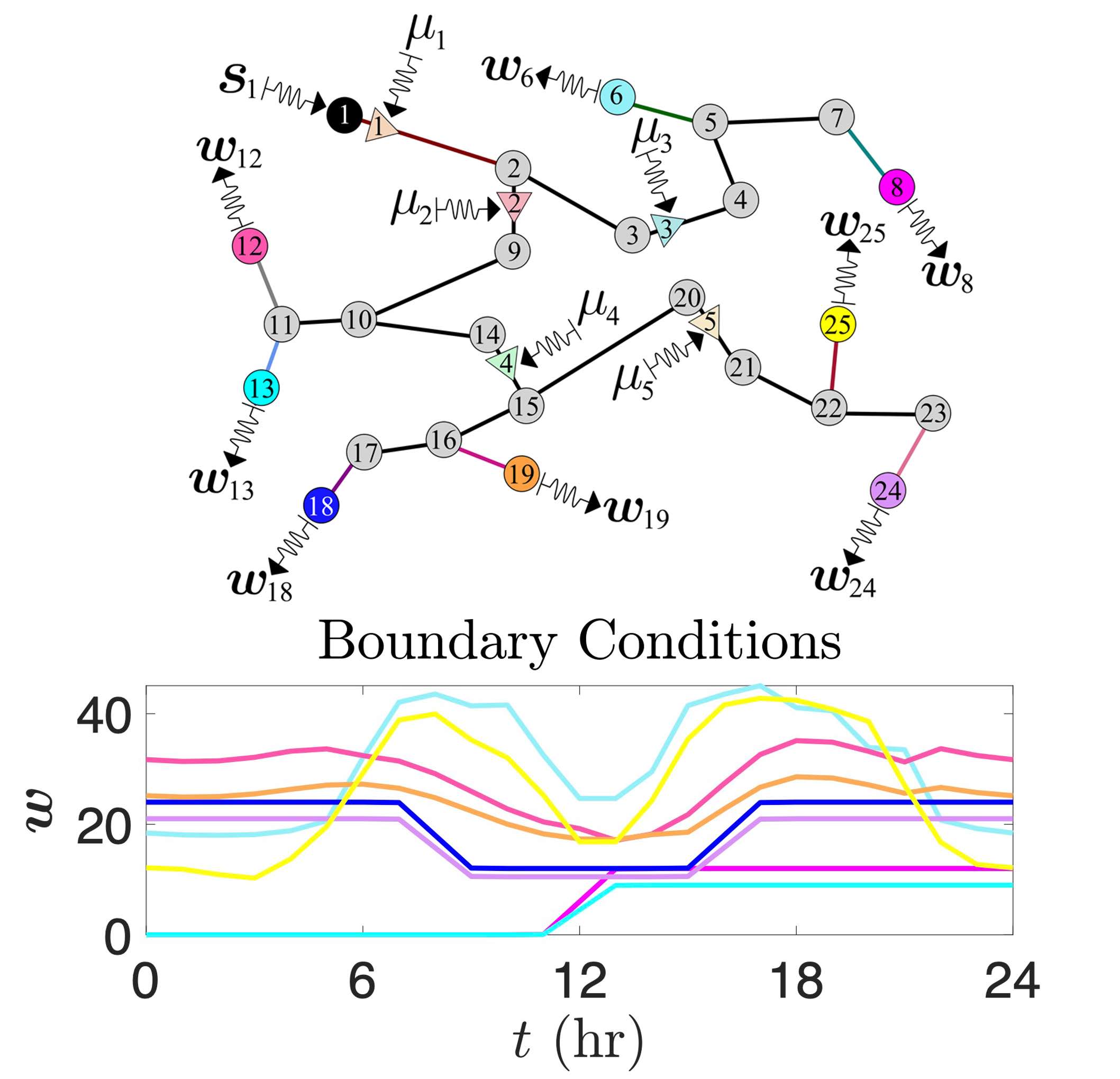}
    \end{subfigure}%
    ~ 
    \begin{subfigure}[t]{0.333\textwidth}
        \centering
        \includegraphics[width=1\linewidth]{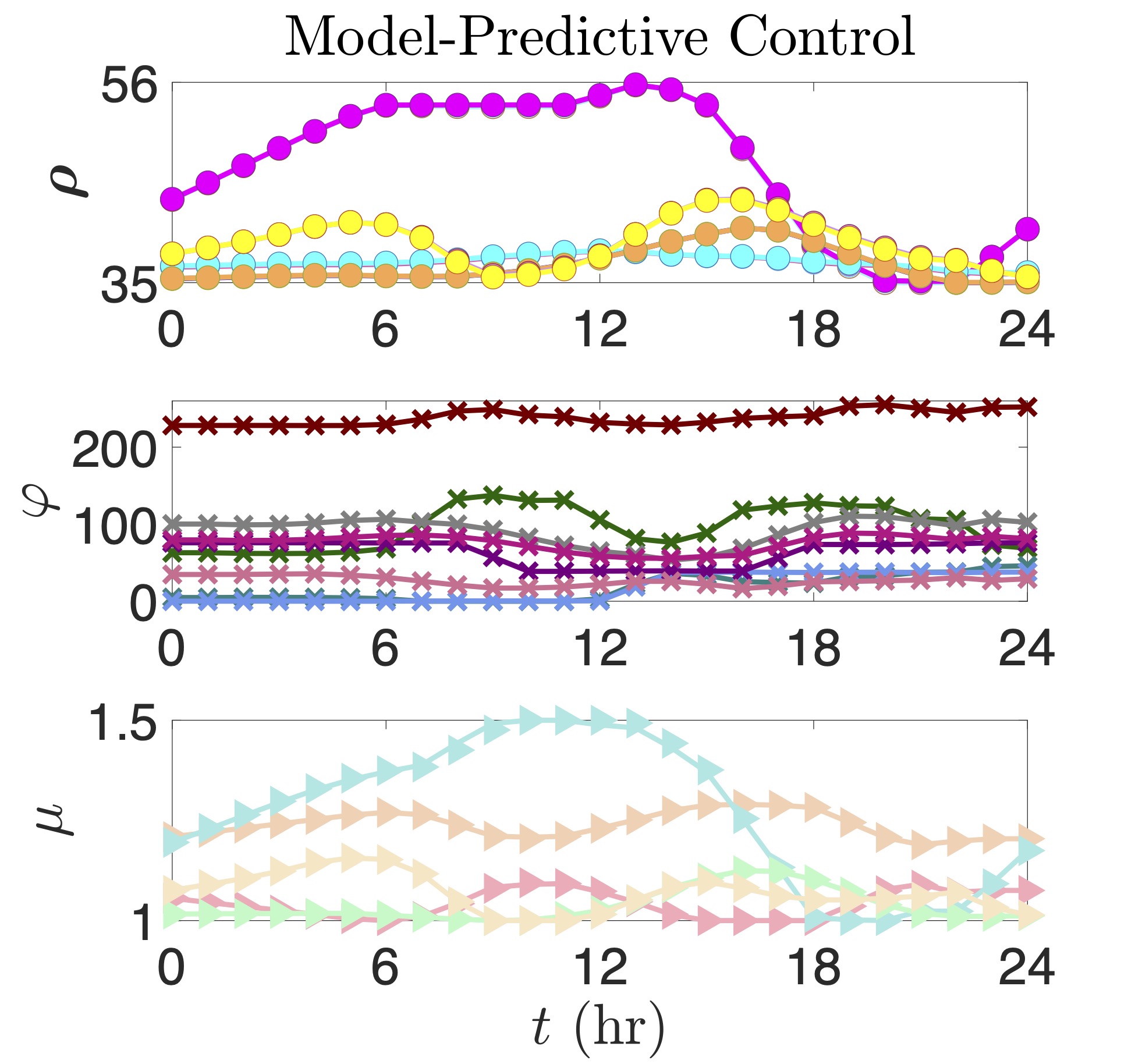}
    \end{subfigure}%
    ~ 
    \begin{subfigure}[t]{0.333\textwidth}
        \centering
        \includegraphics[width=1\linewidth]{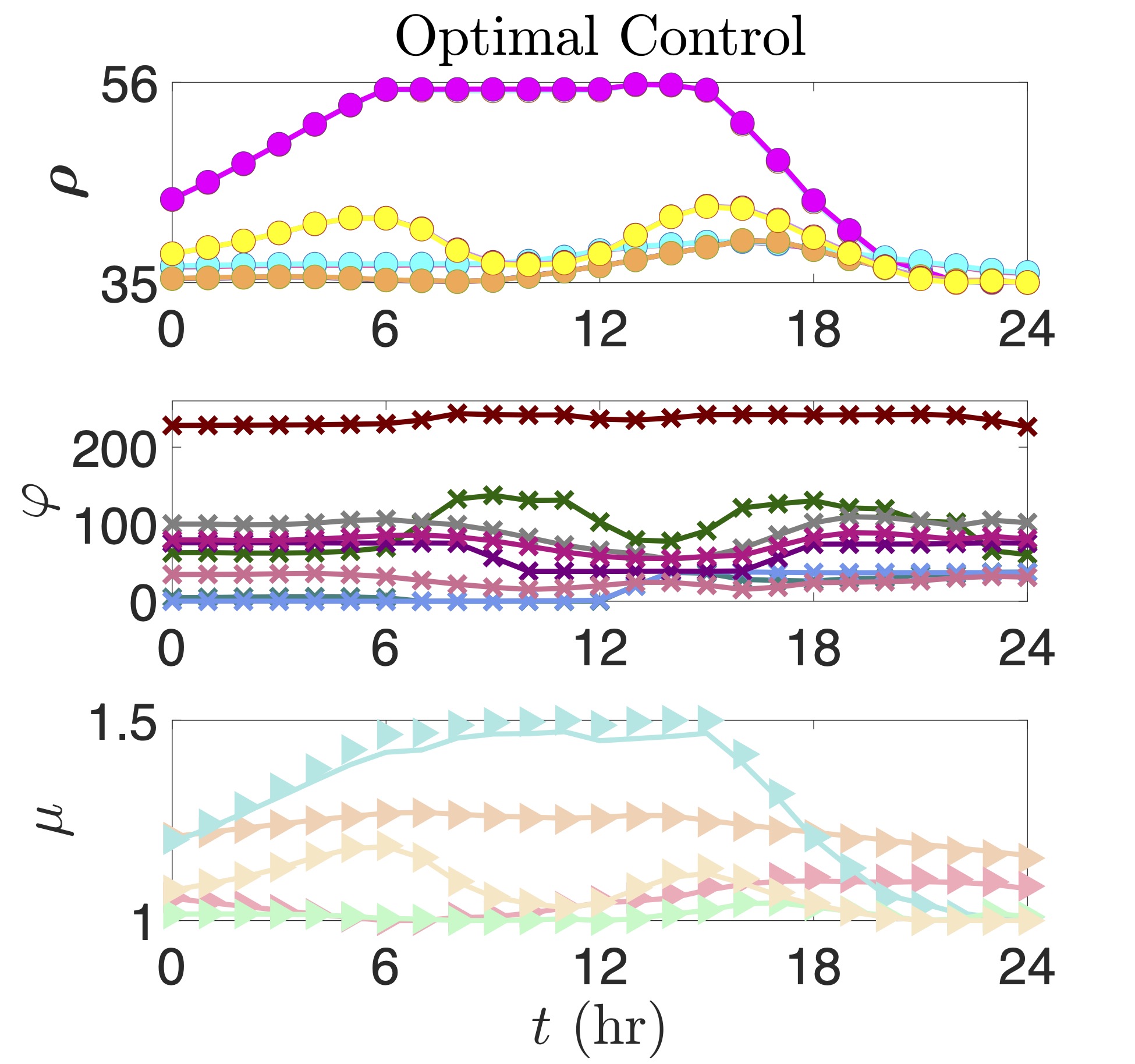}
    \end{subfigure}
    \caption{Left: Network diagram (not to scale) and boundary conditions at color-coordinated nodes \cite{schmidt2017gaslib}.   Center: MPC(60 min) solution. Right: OC(60 min) solution. The sound speed is $\sigma=377.964$ m/s, the total pipe length is 477 km, the friction factor is $\lambda_k=0.01$ for all $k\in \mathcal E$, and the diameters vary between 0.6096 m and 0.9144 m.}
    \label{fig:mpc_25_node}
\end{figure*}

The first network, whose topology is shown in Figure \ref{fig:mpc_cyclic}, was previously used in studies to validate a simulation method based on staggered grid finite difference discretization \cite{GYRYA201934}. The blue node is the sole supply node, with the boundary condition specified as $\bm s_{\text{blue}}=21$ kg/m$^3$. The boundary conditions for loads at withdrawal nodes are depicted in Figure \ref{fig:mpc_cyclic}. The following parameter bounds are used in equation \eqref{eq:ineq}: $\bm \rho^{\min}=21$ kg/m$^3$, $\bm \rho^{\max}=35$ kg/m$^3$, $\varphi^{\min}=0$ kg/m$^2$s, and $\mu_k^{\max}=1.7$ for all $k\in \mathcal C$.  The five edges of the network are discretized into 48 equally space refined segments, so that $\ell_k=5$ km for all $k\in \hat{\mathcal E}$.  This spatial discretization results in a $95\times 95$ state matrix $\overline A$.  

As described above, the four distinct control algorithms will be used to design policies that seek to minimize the energy expended by the three compressor units.  Figure \ref{fig:mpc_cyclic} displays the time evolution of density at the withdrawal nodes, mass flux at the inlets of the edges, and control actions of the compressor units for all four controllers.  Linear and nonlinear solutions are depicted with markers and solid lines, respectively.  Observe that all four controllers determine feasible and qualitatively similar solutions that remain within limitations specified by density, flow, and compressor constraints.  \edit{Furthermore, as documented in Table \ref{tab:cyclic}, the difference in the accumulated compressor energy between MPC and OC is within 5\%.  The energy costs in Table \ref{tab:cyclic} are documented for only the nonlinear controllers because the results for the corresponding linear controllers differ from these by at most 1\%. Table \ref{tab:cyclic} further confirms the error bounds and significant reduction in computational time described by Theorems 3 and 4, respectively. Notably, the computational time for MPC increases linearly with increasing the number of time samples, as expected from the proof of Theorem 4.}

\begin{table}
    \centering
        \caption{\edit{Computational results of MPC and OC for the 5-node cyclic network in Figure \ref{fig:mpc_cyclic}.}}
    \label{tab:cyclic}
    \begin{tabular}{|c|c|c|c|c|c|}
    \hline 
       & Clock & $J$  & $\bm E_{\bm \rho}$\% & $\bm E_{ \varphi}$\% & $\bm E_{\mu}$\% \\
     \hline
     MPC(60\text{min})  & 6.3s &  35.022  & 0.063 & 0.018 & 1.402 \\ 
      \hline
     MPC(30\text{min})  & 9.8s &  35.612  & 0.027 & 0.013 & 0.771 \\ 
           \hline
     MPC(20\text{min})  & 14.9s &  35.188  & 0.011 & 0.012 & 0.478 \\ 
           \hline
     MPC(10\text{min})  & 25.8s &  34.680  & 0.004 & 0.013 & 0.214 \\ 
                \hline
     MPC(5\text{min})  & 42.0s &  34.456  & 0.006 & 0.032 & 0.091 \\ 
           \hline \hline
     OC(60\text{min})  & 117s &  34.085  & 0.677 & 0.078 & 1.529 \\ 
       \hline
     OC(30\text{min})  & 477s &  33.738 & 0.670 & 0.103 & 1.524 \\ 
       \hline
    \end{tabular}
\end{table}

\begin{table}
    \centering
        \caption{\edit{Computational results of MPC and OC for the 25-node tree network in Figure \ref{fig:mpc_25_node}.}}
    \label{tab:25_node}
    \begin{tabular}{|c|c|c|c|c|c|}
    \hline 
       & Clock & $J$  & $\bm E_{\bm \rho}$\% & $\bm E_{ \varphi}$\% & $\bm E_{\mu}$\% \\
     \hline
     MPC(60\text{min})  & 4.9s &  15.655  & 0.078 & 0.151 & 1.215 \\ 
      \hline
     MPC(30\text{min})  & 7.7s &  15.284  & 0.036 & 0.163 & 0.721 \\ 
           \hline
     MPC(20\text{min})  & 11.1s &  15.143  & 0.022 & 0.148 & 0.524 \\ 
           \hline
     MPC(10\text{min})  & 20.3s &  15.057  & 0.013 & 0.211 & 0.349 \\ 
           \hline
    MPC(5\text{min})  & 38.6s &  15.107  & 0.009 & 0.124 & 0.209 \\ 
           \hline \hline
     OC(60\text{min})  & 244s &  15.288  & 0.520 & 0.259 & 2.819 \\ 
       \hline
     OC(30\text{min})  & 1096s &  14.744  & 0.555 & 0.258 & 2.083 \\ 
       \hline
    \end{tabular}
\end{table}

\subsection{25-Node Tree Network}

Our second example involves a 25-node tree network with 5 compressor stations, as depicted in Figure \ref{fig:mpc_25_node}. This network has also been previously studied to demonstrate an OC method similar to the one presented here \cite{zlotnik2015optimal}. Node 1 is the sole supply node, with the boundary condition specified as $\bm s_{1}=35$ kg/m$^3$.  The bounds in equation \eqref{eq:ineq} are defined by $\bm \rho^{\min}=35$ kg/m$^3$, $\bm \rho^{\max}=56$ kg/m$^3$, $\varphi^{\min}=0$ kg/m$^2$s, and $\mu_k^{\max}=1.5$ for all $k\in \mathcal C$.  The 24 edges of the network are discretized into  refined segments in such a way that $\ell_k\le 10$ km for all $k\in \hat{\mathcal E}$.  The size of the resulting state matrix $\overline A$ is $137\times 137$.

Results for MPC and OC are displayed in Figure \ref{fig:mpc_25_node} and are documented in Table \ref{tab:25_node}.  As with the 5-node example, the solutions of this 25-node network exhibit similar qualitative and quantitative behavior among the four controllers.  In particular, the percent difference in accumulated compressor energy between MPC and OC is within 5\%, the errors between linear and nonlinear MPC solutions generally decrease with increasing sampling times, and the computational time for MPC increases linearly with increasing sampling times.  Moreover, the control actions determined by MPC are more responsive to local variations in loads than those determined by OC.  \edit{It is important to note that the solutions presented in Figure \ref{fig:mpc_25_node} for MPC and OC are nearly infeasible in the sense that small increases in subjected loads obstruct the convergence of both programs. Preliminary simulations indicate that the proposed MPC synthesis is more sensitive to such extreme operating conditions and can break down before OC as the conditions approach infeasibility.  These observations lead to two key conclusions.  First, the proposed MPC may not be suitable to predict solutions under extreme operating conditions.  Second, a breakdown of MPC could have significant practical implications for the feasibility of load forecasts and operational safety. Indeed, it may suggest that the physical limitations imposed by the constraints are projected to be severely tested.  Operational limits within which both MPC and OC remain effective are influenced by factors such as network topology, pipeline parameters, and load variations. A detailed analysis of these relationships is beyond the scope of this study.}

%*********************************************************************************************************************************************************************
\section{Conclusions}  \label{sec:conclusion}
%*********************************************************************************************************************************************************************

\edit{As renewable energy continues to drive increased variability in natural gas loads, the ability to efficiently plan for, and adapt to, transient conditions is becoming increasingly important for maintaining system reliability.  This paper has explored the application of adaptive linear systems in a model-predictive control synthesis to optimize transient flows in natural gas transmission pipelines.  Rigorous error bounds derived in this paper offer a framework to assess the validity of modeling natural gas dynamics with adaptive linear systems.  Through extensive numerical simulations, we demonstrate the effectiveness of an adaptive linear MPC formulation to minimize transient compressor energy expenditure while maintaining pipeline flow within operational limits.  One of the key findings is that the proposed MPC offers a significant reduction in computational costs without sacrificing reliability and optimality under usual transient conditions.  Because of the computational efficiency and global optimality provided by convex optimization, this contribution could be significant for establishing baseline conditions for practical use in real-world pipeline operations.

The analysis presented in this paper provides a foundation for further exploration of adaptive linear systems in modeling and controlling transient flows through gas transmission pipelines. While the effectiveness of adaptive linear MPC has been demonstrated on certain systems under sufficiently regular operating conditions, a comprehensive investigation into the feasibility, stability, and robustness of MPC when applied to complex pipeline systems with high load variability and uncertainty is critical. Such an analysis is essential for advancing the integration of MPC into existing pipeline operations.  \revise{Moreover, the error introduced by linearizing the isothermal Euler equations \eqref{eq:gaspde0} remains an open problem that could lead to even more informative guarantees \cite{hall1987linearized, mankbadi1998use, duran2013solution}.} In addition to pure natural gas considerations, it is valuable to note that hydrogen blending has been proposed to reduce the carbon footprint of economies worldwide.  However, optimizing flows of hydrogen blends has been shown to incur considerably higher computational costs compared to homogeneous natural gas  \cite{baker2023optimal, zlotnik2023effects, baker2023transitions, baker2024boundary}. Therefore, pipelines that support clean energy transitions could benefit from an extended MPC specifically tailored to model mixtures of gases. Finally, the online adaptability of MPC could enable more effective communication and dynamic coordination with electric power systems \cite{zlotnik2016coordinated}. These advancements could lead to the development of more efficient and responsive decision-support tools that ensure reliable operation in an increasingly evolving energy landscape.}

\begin{ack}                               % Place acknowledgements
The authors are very grateful to Cody W. Allen, Saif R. Kazi, and Kaarthik Sundar for numerous helpful discussions, and to E. Olga Skowronek for additional discussions and for drawing the network diagrams in Figures \ref{fig:bc_conf}, \ref{fig:mpc_cyclic}, and \ref{fig:mpc_25_node}. This study was supported by the U.S. Department of Energy's Advanced Grid Modeling (AGM) project ``Dynamical Modeling, Estimation, and Optimal Control of Electrical Grid-Natural Gas Transmission Systems,'' the U.S. Department of Energy through the LANL/LDRD Program, and the Center for Non-Linear Studies. Research conducted at Los Alamos National Laboratory is done under the auspices of the National Nuclear Security Administration of the U.S. Department of Energy under Contract No. 89233218CNA000001. L. Baker would like to thank Arizona State University for additional support. S. Shivakumar was supported by the National Science Foundation grants No. 1739990 and 1935453. Report number: LA-UR-23-24798.
\end{ack}

\bibliographystyle{plain}        % Include this if you use bibtex 
\bibliography{autosam}           % and a bib file to produce the 

\begin{thebibliography}{10}

\bibitem{aalto2006real}
Hans Aalto.
\newblock Real-time optimisation of natural gas pipeline systems-a simplified approach.
\newblock {\em IFAC Proceedings Volumes}, 39(14):167--172, 2006.

\bibitem{aalto2015model}
Hans Aalto.
\newblock Model predictive control of natural gas pipeline systems-a case for constrained system identification.
\newblock {\em IFAC-PapersOnLine}, 48(30):197--202, 2015.

\bibitem{abbaspour2008nonisothermal}
Mohammad Abbaspour and Kirby~S. Chapman.
\newblock {Nonisothermal Transient Flow in Natural Gas Pipeline}.
\newblock {\em Journal of Applied Mechanics}, 75(3), 2008.

\bibitem{allen2022gas}
Cody Allen, Roman Zamotorin, Avneet Singh, and Rainer Kurz.
\newblock Gas pipeline transient modeling via transfer functions: Field data validations.
\newblock In {\em PSIG Annual Meeting}, pages PSIG--2205. PSIG, 2022.

\bibitem{baker2021analysis}
Luke Baker, Dieter Armbruster, Anna Scaglione, and Rodrigo~B. Platte.
\newblock Analysis of a model of a natural gas pipeline—a transfer function approach.
\newblock {\em Transactions of Mathematics and Its Applications}, 5(1):tnab002, 2021.

\bibitem{baker2023gas}
Luke~S Baker.
\newblock {\em Gas Mixture Dynamics in Pipeline Networks with a Focus on Linearization and Optimal Control}.
\newblock PhD thesis, Arizona State University, 2023.

\bibitem{baker2023optimal}
Luke~S. Baker, Saif~R. Kazi, Rodrigo~B. Platte, and Anatoly Zlotnik.
\newblock Optimal control of transient flows in pipeline networks with heterogeneous mixtures of hydrogen and natural gas.
\newblock In {\em 2023 American Control Conference (ACC)}, pages 1221--1228. IEEE, 2023.

\bibitem{baker2023transitions}
Luke~S. Baker, Saif~R. Kazi, and Anatoly Zlotnik.
\newblock Transitions from monotonicity to chaos in gas mixture dynamics in pipeline networks.
\newblock {\em PRX Energy}, 2(3):033008, 2023.

\bibitem{baker2024boundary}
Luke~S. Baker and Anatoly Zlotnik.
\newblock Boundary control for suppressing chaotic response to dynamic hydrogen blending in a gas pipeline.
\newblock {\em IFAC-PapersOnLine}, 58(5):72--77, 2024.

\bibitem{behbahani2010accuracy}
Morteza Behbahani-Nejad and Younes Shekari.
\newblock The accuracy and efficiency of a reduced-order model for transient flow analysis in gas pipelines.
\newblock {\em Journal of Petroleum Science and Engineering}, 73(1-2):13--19, 2010.

\bibitem{beylin2020fast}
Aleksandr Beylin, Aleksandr~M. Rudkevich, and Anatoly Zlotnik.
\newblock Fast transient optimization of gas pipelines by analytic transformation to linear programs.
\newblock In {\em PSIG Annual Meeting}, pages PSIG--2003, 2020.

\bibitem{borraz2013optimization}
Conrado Borraz-S{\'a}nchez and Dag Haugland.
\newblock Optimization methods for pipeline transportation of natural gas with variable specific gravity and compressibility.
\newblock {\em Top}, 21:524--541, 2013.

\bibitem{de2000gas}
Daniel De~Wolf and Yves Smeers.
\newblock The gas transmission problem solved by an extension of the simplex algorithm.
\newblock {\em Management Science}, 46(11):1454--1465, 2000.

\bibitem{domschke2011combination}
Pia Domschke, Bjorn Gei{\ss}ler, Oliver Kolb, Jens Lang, Alexander Martin, and Antonio Morsi.
\newblock Combination of nonlinear and linear optimization of transient gas networks.
\newblock {\em INFORMS Journal on Computing}, 23(4):605--617, 2011.

\bibitem{duran2013solution}
Ignacio Duran and Stephane Moreau.
\newblock Solution of the quasi-one-dimensional linearized euler equations using flow invariants and the magnus expansion.
\newblock {\em J. Fluid Mechanics}, 723:190--231, 2013.

\bibitem{ehrhardt2005nonlinear}
Klaus Ehrhardt and Marc~C. Steinbach.
\newblock Nonlinear optimization in gas networks.
\newblock In {\em Modeling, simulation and optimization of complex processes}, pages 139--148. Springer, 2005.

\bibitem{feldman1988optimization}
Michael~Bos Feldman.
\newblock Optimization of gas transmission systems using linear programming.
\newblock In {\em PSIG Annual Meeting}, pages PSIG--8809. PSIG, 1988.

\bibitem{gopalakrishnan2013economic}
Ajit Gopalakrishnan and Lorenz~T. Biegler.
\newblock Economic nonlinear model predictive control for periodic optimal operation of gas pipeline networks.
\newblock {\em Computers \& Chemical Engineering}, 52:90--99, 2013.

\bibitem{grundel2013computing}
Sara Grundel, Nils Hornung, Bernhard Klaassen, Peter Benner, and Tanja Clees.
\newblock Computing surrogates for gas network simulation using model order reduction.
\newblock In {\em Surrogate-Based Modeling and Optimization}, pages 189--212. Springer, 2013.

\bibitem{gugat2020closed}
Martin Gugat, Falk~M. Hante, and Li~Jin.
\newblock Closed loop control of gas flow in a pipe: stability for a transient model.
\newblock {\em at-Automatisierungstechnik}, 68(12):1001--1010, 2020.

\bibitem{gugat2012well}
Martin Gugat, Michael Herty, Axel Klar, G{\"u}nter Leugering, and Veronika Schleper.
\newblock Well-posedness of networked hyperbolic systems of balance laws.
\newblock {\em Constrained Optimization and Optimal Control for Partial Differential Equations}, 160:123--146, 2012.

\bibitem{gugat2016stationary}
Martin Gugat, G{\"u}nter Leugering, and Falk Hante.
\newblock Stationary states in gas networks.
\newblock {\em Networks and Heterogeneous Media}, 10(2):295--320, 2016.

\bibitem{gugat2018mip}
Martin Gugat, G{\"u}nter Leugering, Alexander Martin, Martin Schmidt, Mathias Sirvent, and David Wintergerst.
\newblock Mip-based instantaneous control of mixed-integer pde-constrained gas transport problems.
\newblock {\em Computational Optimization and Applications}, 70:267--294, 2018.

\bibitem{gugat2017isothermal}
Martin Gugat and Stefan Ulbrich.
\newblock The isothermal {Euler} equations for ideal gas with source term: Product solutions, flow reversal and no blow up.
\newblock {\em Journal of Mathematical Analysis and Applications}, 454(1):439--452, 2017.

\bibitem{GYRYA201934}
Vitaliy Gyrya and Anatoly Zlotnik.
\newblock An explicit staggered-grid method for numerical simulation of large-scale natural gas pipeline networks.
\newblock {\em Applied Mathematical Modelling}, 65:34--51, 2019.

\bibitem{hall1987linearized}
Kenneth~Charles Hall.
\newblock {\em A linearized Euler analysis of unsteady flows in turbomachinery}.
\newblock PhD thesis, Massachusetts Institute of Technology, 1987.

\bibitem{himpe2021model}
Christian Himpe, Sara Grundel, and Peter Benner.
\newblock Model order reduction for gas and energy networks.
\newblock {\em Journal of Mathematics in Industry}, 11(1):1--46, 2021.

\bibitem{horn2012matrix}
Roger~A. Horn and Charles~R. Johnson.
\newblock {\em Matrix Analysis}.
\newblock Cambridge University Press, 2012.

\bibitem{khalil2002nonlinear}
Hassan~K. Khalil.
\newblock {\em Nonlinear Systems}.
\newblock Pearson Education. Prentice Hall, 2002.

\bibitem{kralik1984modeling}
Jaroslav Kr{\'a}lik, Petr Stiegler, Zden{\v{e}}k Vostr{\'u}, and Ji{\v{r}}i Z{\'a}vorka.
\newblock Modeling the dynamics of flow in gas pipelines.
\newblock {\em IEEE Transactions on Systems, Man, and Cybernetics}, (4):586--596, 1984.

\bibitem{kralik1984universal}
Jaroslav Kr{\'a}lik, Petr Stiegler, Zden{\v{e}}k Vostr{\`y}, and Ji{\v{r}}{\'\i} Z{\'a}vorka.
\newblock A universal dynamic simulation model of gas pipeline networks.
\newblock {\em IEEE Transactions on Systems, Man, and Cybernetics}, (4):597--606, 1984.

\bibitem{liu2011coordinated}
Cong Liu, Mohammad Shahidehpour, and Jianhui Wang.
\newblock Coordinated scheduling of electricity and natural gas infrastructures with a transient model for natural gas flow.
\newblock {\em Chaos: An Interdisciplinary Journal of Nonlinear Science}, 21(2):025102, 2011.

\bibitem{luongo1986efficient}
Cesar~A. Luongo.
\newblock An efficient program for transient flow simulation in natural gas pipelines.
\newblock In {\em {PSIG Annual Meeting}}, pages PSIG--8605, 1986.

\bibitem{mankbadi1998use}
R.~R. Mankbadi, R.~Hixon, S.-H. Shih, and L.~A. Povinelli.
\newblock Use of linearized euler equations for supersonic jet noise prediction.
\newblock {\em AIAA journal}, 36(2):140--147, 1998.

\bibitem{ivan2005}
Ivan Marić, Antun Galovic, and Tomislav Smuc.
\newblock Calculation of natural gas isentropic exponent.
\newblock {\em Flow Measurement and Instrumentation}, 16:13--20, 03 2005.

\bibitem{marques1988line}
Dardo Marqu{\'e}s and Manfred Morari.
\newblock On-line optimization of gas pipeline networks.
\newblock {\em Automatica}, 24(4):455--469, 1988.

\bibitem{mayne2000constrained}
David~Q. Mayne, James~B. Rawlings, Christopher~V. Rao, and Pierre O.~M. Scokaert.
\newblock Constrained model predictive control: Stability and optimality.
\newblock {\em Automatica}, 36(6):789--814, 2000.

\bibitem{misra2014optimal}
Sidhant Misra, Michael~W. Fisher, Scott Backhaus, Russell Bent, Michael Chertkov, and Feng Pan.
\newblock Optimal compression in natural gas networks: A geometric programming approach.
\newblock {\em {Transactions on Control of Network Systems}}, 2(1):47--56, 2014.

\bibitem{osiadacz1984simulation}
Andrzej Osiadacz.
\newblock Simulation of transient gas flows in networks.
\newblock {\em {International Journal for Numerical Methods in Fluids}}, 4(1):13--24, 1984.

\bibitem{osiadacz2001comparison}
Andrzej~J. Osiadacz and Maciej Chaczykowski.
\newblock Comparison of isothermal and non-isothermal pipeline gas flow models.
\newblock {\em Chemical Engineering Journal}, 81(1-3):41--51, 2001.

\bibitem{percell1987steady}
Peter~B. Percell and Michael~J. Ryan.
\newblock Steady state optimization of gas pipeline network operation.
\newblock In {\em PSIG Annual Meeting}, pages PSIG--8703, 1987.

\bibitem{Yue20}
Yue Qiu, Sara Grundel, Martin Stoll, and Peter Benner.
\newblock Efficient numerical methods for gas network modeling and simulation.
\newblock {\em Networks and Heterogeneous Media}, 15(4):653--679, 2020.

\bibitem{schmidt2017gaslib}
Martin Schmidt, Denis A{\ss}mann, Robert Burlacu, Jesco Humpola, Imke Joormann, Nikolaos Kanelakis, Thorsten Koch, et~al.
\newblock Gaslib—a library of gas network instances.
\newblock {\em Data}, 2(4):40, 2017.

\bibitem{srinivasan2022numerical}
Shriram Srinivasan, Kaarthik Sundar, Vitaliy Gyrya, and Anatoly Zlotnik.
\newblock Numerical solution of the steady-state network flow equations for a non-ideal gas.
\newblock {\em IEEE Transactions on Control of Network Systems}, 2022.

\bibitem{steinbach2007pde}
Marc~C. Steinbach.
\newblock On {PDE} solution in transient optimization of gas networks.
\newblock {\em Journal of Computational and Applied Mathematics}, 203(2):345--361, 2007.

\bibitem{sundar2018state}
Kaarthik Sundar and Anatoly Zlotnik.
\newblock State and parameter estimation for natural gas pipeline networks using transient state data.
\newblock {\em IEEE Transactions on Control Systems Technology}, 27(5):2110--2124, 2018.

\bibitem{thorley1987unsteady}
A.~R.~D. Thorley and C.~H. Tiley.
\newblock Unsteady and transient flow of compressible fluids in pipelines—a review of theoretical and some experimental studies.
\newblock {\em International Journal of Heat and Fluid Flow}, 8(1):3--15, 1987.

\bibitem{thulasiraman2011graphs}
Krishnaiyan Thulasiraman and Madisetti N.~S. Swamy.
\newblock {\em Graphs: theory and algorithms}.
\newblock John Wiley \& Sons, 2011.

\bibitem{vaidya1989speeding}
Pravin~M. Vaidya.
\newblock Speeding-up linear programming using fast matrix multiplication.
\newblock In {\em 30th annual symposium on foundations of computer science}, pages 332--337. IEEE Computer Society, 1989.

\bibitem{wong1968optimization}
P.~L. Wong and R.~Larson.
\newblock Optimization of natural-gas pipeline systems via dynamic programming.
\newblock {\em {IEEE Transactions on Automatic Control}}, 13(5):475--481, 1968.

\bibitem{wu2000model}
S.~Wu, R.~Z. Rios-Mercado, E.~A. Boyd, and L.~R. Scott.
\newblock Model relaxations for the fuel cost minimization of steady-state gas pipeline networks.
\newblock {\em {Mathematical and Computer Modelling}}, 31(2-3):197--220, 2000.

\bibitem{zecchin}
Aaron~C. Zecchin, Angus~R. Simpson, Martin~F. Lambert, Langford~B. White, and John~P. V{\'\i}tkovsk{\`y}.
\newblock Transient modeling of arbitrary pipe networks by a laplace-domain admittance matrix.
\newblock {\em J. Engineering Mechanics}, 135(6):538--547, 2009.

\bibitem{zhang2020transient}
Shixuan Zhang, Sheng Liu, Tianhu Deng, and Zuo-Jun~Max Shen.
\newblock Transient-state natural gas transmission in gunbarrel pipeline networks.
\newblock {\em INFORMS Journal on Computing}, 32(3):697--713, 2020.

\bibitem{grail2018data}
Anatoly Zlotnik.
\newblock Gas reliability analysis integrated library, 2018.
\newblock \verb"https://github.com/lanl-ansi/grail/".

\bibitem{zlotnik2015optimal}
Anatoly Zlotnik, Michael Chertkov, and Scott Backhaus.
\newblock Optimal control of transient flow in natural gas networks.
\newblock In {\em Conference on Decision and Control}, pages 4563--4570. IEEE, 2015.

\bibitem{zlotnik2015model}
Anatoly Zlotnik, Sergey Dyachenko, Scott Backhaus, and Michael Chertkov.
\newblock Model reduction and optimization of natural gas pipeline dynamics.
\newblock In {\em Dynamic Systems and Control Conference}, volume 57267, 2015.

\bibitem{zlotnik2023effects}
Anatoly Zlotnik, Saif~R. Kazi, Kaarthik Sundar, Vitaliy Gyrya, Luke Baker, Mo~Sodwatana, and Yan Brodskyi.
\newblock Effects of hydrogen blending on natural gas pipeline transients, capacity, and economics.
\newblock In {\em PSIG Annual Meeting}, pages PSIG--2312, 2023.

\bibitem{zlotnik2016coordinated}
Anatoly Zlotnik, Line Roald, Scott Backhaus, Michael Chertkov, and G{\"o}ran Andersson.
\newblock Coordinated scheduling for interdependent electric power and natural gas infrastructures.
\newblock {\em IEEE Transactions on Power Systems}, 32(1):600--610, 2016.

\bibitem{zlotnik2017economic}
Anatoly Zlotnik, Aleksandr~M. Rudkevich, Evgeniy Goldis, Pablo~A. Ruiz, Michael Caramanis, Richard Carter, Scott Backhaus, Richard Tabors, Richard Hornby, and Daniel Baldwin.
\newblock Economic optimization of intra-day gas pipeline flow schedules using transient flow models.
\newblock In {\em PSIG Annual Meeting}, pages PSIG--1715, 2017.

\end{thebibliography}
                                 % bibliography (preferred). The
                                 % correct style is generated by
                                 % Elsevier at the time of printing.

%\begin{thebibliography}{99}     % Otherwise use the  
                                 % thebibliography environment.
                                 % Insert the full references here.
                                 % See a recent issue of Automatica 
                                 % for the style.
%  \bibitem[Heritage, 1992]{Heritage:92}
%     (1992) {\it The American Heritage. 
%     Dictionary of the American Language.}
%     Houghton Mifflin Company.
%  \bibitem[Able, 1956]{Abl:56}
%     B.~C.~Able (1956). Nucleic acid content of macroscope. 
%     {\it Nature 2}, 7--9. 
%  \bibitem[Able {\em et al.}, 1954]{AbTaRu:54}   
%     B.~C. Able, R.~A. Tagg, and M.~Rush (1954).
%     Enzyme-catalyzed cellular transanimations.
%     In A.~F.~Round, editor, 
%     {\it Advances in Enzymology Vol. 2} (125--247). 
%     New York, Academic Press.
%  \bibitem[R.~Keohane, 1958]{Keo:58}
%     R.~Keohane (1958).
%     {\it Power and Interdependence: 
%     World Politics in Transition.}
%     Boston, Little, Brown \& Co.
%  \bibitem[Powers, 1985]{Pow:85}
%     T.~Powers (1985).
%     Is there a way out?
%     {\it Harpers, June 1985}, 35--47.

%\end{thebibliography}

% \appendix
% \section{A summary of Latin grammar}    % Each appendix must have a short title.
% \section{Some Latin vocabulary}         % Sections and subsections are supported  
%                                         % in the appendices.
\end{document}